\RequirePackage{fix-cm}
\documentclass[smallcondensed]{svjour3}     % onecolumn (ditto)
\smartqed  % flush right qed marks, e.g. at end of proof
\usepackage{graphicx}
 \usepackage{mathptmx}      % use Times fonts if available on your TeX system
\usepackage{amsmath}
\usepackage{amssymb}
\usepackage{color}
\usepackage{empheq}

\def \n {{\text{n}}}

\def \Qs {{B}_{\varepsilon}}
\def \Qc {{D}_{\varepsilon}}
\def \Qcc {{A}_{\varepsilon}}
\def \ve{\varepsilon}
\def \be{\begin{equation}}
\def \ee{\end{equation}}

\newtheorem{corrolary}{Corrolary}

\begin{document}

\title{Phase field model of cell motility: sharp interface limit in sub-critical case
\thanks{The work of LB was supported by NSF grant DMS-1106666.
The work of VR and MP was partially supported by NSF grant DMS-1106666. 
The authors are grateful to I. Aronson and F. Ziebert for useful discussions on  the phase field model of cell motility introduced in their paper.}
}


\author{Leonid Berlyand        \and
        Volodymyr Rybalko  \and Mykhailo Potomkin%etc.
}


\institute{L. Berlyand \at
             Department of Mathematics, 
	The Pennsylvania State University, 
	University Park, PA 16802, USA\\    
              \email{berlyand@math.psu.edu}           %  \\
           \and
           V. Rybalko \at
              Mathematical Division, B. Verkin Institute for Low Temperature, Physics and Engineering
	of National Academy of Sciences of Ukraine,
	47 Lenin Ave., 61103 Kharkiv, Ukraine\\
	\email{vrybalko@ilt.kharkov.ua}
	\and 
	M. Potomkin \at 
	The Pennsylvania State University, 
	University Park, PA 16802, USA,\\
	\email{mup20@ucs.psu.edu}   
}

\date{}

\maketitle

\begin{abstract}

We consider a system of two PDEs  arising in modeling of motility of eukariotic cells on substrates.  This system consists of  the Allen-Cahn equation for the  scalar phase field   function coupled with another vectorial  parabolic equation for the orientation  of the actin filament network.

 The two key  properties  of this system are   (i) presence of gradients in the coupling terms (gradient coupling) and  (ii) mass (volume) preservation constraints.
 We  first prove that the sharp interface property  of initial conditions is preserved in time. Next we formally derive the equation of the motion of the interface, which  is  the  mean curvature  motion perturbed by a nonlinear  term  that appears due  to the properties (i)-(ii).  This  novel term leads to surprising  features of the the motion of the interface. 
 
 Because of these properties  maximum principle and classical  comparison techniques   do not  apply to this system.  Furthermore, the system can not be written in  a form of gradient flow, which is why recently developed  $\Gamma$-convergence techniques also  can not be used for the justification of the formal derivation.   Such justification is presented in a one-dimensional model problem and it leads to a stability result in  a  class of ``sharp interface" initial data. 
 
\end{abstract}

\section{Introduction}
\label{intro}

The problem of cell motility  has been a classical subject in biology for  several centuries. It  dates back to the celebrated discovery by  van Leeuwenhoek in 17th century  who  drastically improved microscope   to the extent that he was able   to  observe motion of single celled organisms that moved due contraction and extension.  Three centuries later this problem  continues to  attract  the  attention  of  biologists,   biophysicists  and, more recently,  applied mathematicians.  A comprehensive  review  of  the mathematical modeling of cell motility can be found in \cite{Mog09}.

This work is motivated by the  problem of  motility (crawling motion) of eukariotic cells on substrates.  The network of  actin  (protein) filaments (which is  a part of cytoskeleton  in such cells)   plays an important role in cell motility.   We are concerned with  the cell motility caused by extension of front of the cell due to polymerization of the actin filaments  and contraction of the back of the cell due to detachment of  these filaments.  Modeling of  this process in full generality is at present a formidable challenge because  several  important  biological   ingredients  (e.g., regulatory pathways \cite{Mog09}) are not yet  well understood.  

However, in recent biophysical studies  several simplified  {\it phase field models}  have been proposed. Simulations performed for these models demonstrated   good agreement  with experiments (e.g., \cite{ShaRapLev10,ZieSwaAra11} and references their in).  Recall that phase field models are typically used to describe the  evolution of an interface between two phases (e.g., solidification or viscous fingering). The key ingredient of such models is  an auxiliary  scalar field,  which takes two  different values in   domains  describing  the two phases (e.g., $1$ and $0$) with a diffuse interface of a small  (non zero) width.  The advantage of such an approach is that it allows us to consider one set of PDEs in the whole domain  occupied by both phases and therefore avoids an issue of coupling two different  sets of PDEs in each phase, which is typically quite  complicated in both simulations and analysis.  

In this work we present rigorous mathematical analysis of the  2D phase field model proposed in \cite{ZieSwaAra11} that consists of a system of  two PDEs for the  phase field function  and  orientation vector with an integral mass conservation constraint.  This model  can be rewritten in a simplified  form suitable  for  asymptotical analysis,  so that all key features of the qualitative behavior are preserved, which can be seen from a comparison of simulations  from \cite{ZieSwaAra11} with our analytical results.  First, in \cite{ZieSwaAra11}  the integral mass conservation constraint is introduced in  the PDE system by adding  a penalization parameter  into  the double-well  potential  (formulas (2.2), (2.5)-(2.6) in  \cite{ZieSwaAra11}).  We introduce this constraint via a dynamic Lagrange multiplier $\lambda(t)$ defined below, which  provides the same qualitative  behavior of solutions. Second, for technical simplicity we drop two terms in the second  equation (for polarization) in the phase field system  \cite{ZieSwaAra11}, since our analysis shows that these terms can be incorporated  with minor changes in both the results  and the techniques. Thirdly,  in order to study the long term behavior of the system,  we perform   the diffusive scaling ($t\mapsto \ve^2 t$, $x \mapsto \ve x$). Indeed, the crawling motion is very slow and  time  variable needs to be ``accelerated". Thus, we arrived at  the following  system  of parabolic PDEs  for a  scalar  phase field function 
$\rho_\ve$ and the orientation vector $P_\ve$: 
\begin{equation}
\frac{\partial \rho_\ve}{\partial t}=\Delta \rho_\ve
-\frac{1}{\ve^2}W^{\prime}(\rho_\ve)%\rho(\rho-1)(\rho-1/2)
-P_\ve\cdot \nabla \rho_\ve +\lambda_\ve(t)
  \text{ in }\
\Omega, 
\label{eq1}
\end{equation}
\begin{equation}
\frac{\partial P_\ve}{\partial t}=\ve\Delta P_\ve -\frac{1}{\ve}P_\ve
 -\beta \nabla \rho_\ve
 %-\textcolor{red}{\gamma_\ve} (\nabla \rho\cdot P)\, P
 \qquad\text{in}\ \Omega.
\label{eq2}
\end{equation}
On the boundary $\partial \Omega$ we impose the Neumann  and the Dirichlet boundary conditions respectively $\partial_\nu \rho_\ve=0$ and   $P_\ve=0$, 
where $\Omega \subset \mathbb{R}^2$ is a bounded smooth domain.

Equation \eqref{eq1} is a  perturbation of the following Allen-Cahn equation \cite{AllCah97}:
\begin{equation} \label{AC}
\frac{\partial \rho_\ve}{\partial t}=\Delta \rho_\ve
-{\frac{1}{\ve^2}}W^\prime(\rho_\ve).%(1-\rho)(1/2-\rho)\rho.
\end{equation}
 The latter equation  is a scalar version  of  the  celebrated Ginzburg-Landau equation and it plays a fundamental role in  mathematical modeling of phase transitions. It consists of  the standard linear  parabolic equation and a nonlinear lower order  term, which is the derivative of a smooth double equal  well potential   
  \begin{equation}\label{doublewell}
  W(\rho)=\frac{1}{4} \rho^2 (1-\rho)^2.
  \end{equation}
    Equation \eqref{AC}  was introduced to  model  the motion of phase-antiphase boundary (interface) between two grains in a solid material.
%JSP, 1994, v.77, Fife and Lacey.  
%the motion of a boundary that separates two phases of a polycrystalline materials. 
    Analysis of \eqref{AC} as $\ve \to 0$  led to the asymptotic solution that takes values $0$ and $1$ in the domains  corresponding to  two phases separated by an interface of the  width  of order $\ve$,  the so-called  {\it sharp interface}.  Furthermore, it was shown that this sharp interface  exhibits  the  {\it mean curvature motion}.   Recall that  in this motion the normal component of the velocity of  each point of the surface  is equal to the mean curvature  of the surface at this point. This  motion  has been  extensively studied  in the geometrical community (e.g., \cite{Ham82,Hui84,Gra87,Bra78} and references therein).  It also  received significant attention  in
  PDE  literature.  Specifically \cite{CheGigGot91}  and   \cite{EvaSpr91}  established existence of global viscosity  solutions (weak solutions) for the  mean curvature flow.  The mean curvature motion of the interface in the limit $\ve \to 0$  was formally derived in \cite{RubSterKel89}, \cite{Fif88} and  then justified in \cite{EvaSonSou92}  by using the viscosity  solutions techniques.

  % \cite{EvansSoner Souganidis, Evans Spruck, Sternberg-Rubinstein curve shortening,}  communities. 

 Note that equation \eqref{AC} is closely related to another well-known model of phase separation, the so-called Cahn-Hilliard equation \cite{CahHil58}, which  is a forth order  reaction diffusion equation that models how two components of a binary fluid spontaneously separate and form  domains  of  two pure  fluids.

There are    two   distinguishing  features  in the  problem  \eqref{eq1}-\eqref{eq2}: coupling and 
a  nonlocal mass conservation constraint.
%  a special choice of parameters in singular perturbations of both equations.
We first comment on the  coupling. Note that  another prominent biological  FitzHugh-Nagumo  
model  has similar coupling feature but in \eqref{eq1}-\eqref{eq2} the {\it coupling occurs via spatial gradients (gradient  coupling)}  of the unknown functions   where as in FitzHugh-Nagumo \cite{RocGeoGiu00}, \cite{SorSou96} the two equations are coupled via lower order terms (unknown functions  rather than their derivatives).  There are several  phase field models for  Allen-Cahn  (also Cahn-Hilliard) equation   coupled with another parabolic equation via lower order terms \cite{EvaSonSou92,Che94}. Our analysis shows that  the gradient  coupling results in novel mathematical features such as the following  nonlinear nonlocal equation for the  velocity of the interface curve $\Gamma(t)$ derived below: 
\begin{equation}\label{motion1}
V=\kappa +\frac{\beta}{c_0}\Phi(V)-\frac{1}{|\Gamma(t)|}\int_{\Gamma(t)}\left(\kappa
+\frac{\beta}{c_0}\Phi(V)\right)\,  ds. 
\end{equation}  
 Here $V$ stands for the normal velocity of $\Gamma(t)$ with respect to the  inward normal, and $\kappa$ 
denotes the curvature of $\Gamma(t)$,   $c_0$ is a constant determined by the potential $W$ ($c_0=\sqrt{3/2}$ for the specific choice \eqref{doublewell} of the double-well potential), $|\Gamma(t)|$ is the  curve length,  and function $\Phi(V)$  is  given by \eqref{def_for_motion1}.
% can be  computed explicitly.   

Next  note  that the term $\lambda_\ve(t)$ in \eqref{eq1} is a Lagrange multiplier responsible for the  
{\it volume constraint} (conservation mass in the original physical problem \cite{ZieSwaAra11})  and it   has the following form

\begin{equation} \label{lagrange}
\lambda_\ve(t)=\frac{1}{|\Omega|}\int_\Omega\left(\frac{1}{\ve^2}W^\prime(\rho_\ve)
+ P_\ve\cdot \nabla \rho_\ve \right)\, dx
\end{equation}

Solutions of stationary Allen-Cahn equation with the volume constraint were studied in \cite{Mod86} by $\Gamma$-convergence  techniques applied to the stationary variational   problem corresponding  to\eqref{AC}. It was established that the $\Gamma$-limiting functional is the interface perimeter  (curve length in 2D  or surface area  in higher dimensions).  Subsequently in the work 
\cite{RubSte92}   an evolutionary  reaction-diffusion  equation with double-well potential   and  nonlocal term that describes  the volume constraint was studied.  The following asymptotic formula   for evolution of interface  $\Gamma$ in the form of volume preserving mean curvature flow was formally derived in \cite{RubSte92}
\begin{equation} \label{rs_1992}
V=\kappa -\frac{1}{|\Gamma(t)|}\int_{\Gamma(t)}\kappa\, ds
%\qquad\text{Note that}\ V\
%\text{is the velocity w.r.t. inward pointing normal}
\end{equation}
%Here $V$  is the normal velocity of the interface, $\kappa$ is  its mean curvature, and 
 Formula \eqref{rs_1992} was rigorously justified in the radially symmetric  case in \cite{BroSto97} and in general case in \cite{CheHilLog10}.

There are three main approaches to the study of asymptotic behavior (sharp interface limit) of solutions of phase field equations and systems. 

When comparison principle for solutions applies,  a PDE approach based on viscosity solutions techniques  was successfully used in \cite{EvaSonSou92,BarSonSou93,Sou02,LioKimSle04,BarLio03,Gol97}.
% \cite{ A Geometrical Approach to Front Propagation Problems in Bounded Domains with Neumann-Type Boundary Conditions Guy Barles and Francesca Da Lio,2004}.
 This approach can not be applied to  the system \eqref{eq1}-\eqref{eq2}, because of the gradient coupling and nonlocal multiplier $\lambda_\ve(t)$. It is an open issue to  introduce weak (e.g., viscosity or Brakke type) solutions in problems with constraints.   Furthermore,  since there is no \textcolor{black}{comparison  principle}, the only technique available for  the justification of the sharp interface limit  is energy bounds  which become quite difficult due to the  coupling and the volume preservation.
 
  Another technique  used in such problems  is $\Gamma-$convergence (\textcolor{black}{see \cite{Ser10} and references therein}). It  also does not work for the system \eqref{eq1}-\eqref{eq2}. Standard Allen-Cahn equation \eqref{AC} is a gradient flow  (in $L^2$ metric) with GL energy functional, which is why one can use     
 the $\Gamma-$convergence approach. However, there is no energy functional such that problem \eqref{eq1}-\eqref{eq2}   can be  written as a  gradient flow.

% Solutions of stationary Alan-Cahn equation with the volume constraint were studied in \cite{Mod86} by $\Gamma$-convergence  techniques applied to the corresponding variational problem. It was established that the $\Gamma$-limit minimizes the surface area
%
% Bronsard  volume conserved but  evolution problem corresponds  approximately to $\mu=0$ but $\delta$ still present.
%Sternberg and others: evolution but no constrain

 As explained above  the gradient coupling and  the  volume constraint are the  key features of the problem  \eqref{eq1}, \eqref{eq2}, \eqref{lagrange} and they  led   to  both novel results and analysis techniques.   Specifically, the  objectives of  our study  are three fold:

 (i) To show that there is no  finite time blow up and the sharp interface property of the initial data propagates in time.
 
 (ii)  To investigate how the  gradient coupling   combined with the nonlocal  volume constraints  changes the  limiting equation of the interface motion.
 
 (iii) To develop novel techniques for the justification of the limiting  equation  of the interface motion that, in particular, includes  rigorous derivation of asymptotic  expansion for the solution of the problem \eqref{eq1}-\eqref{eq2}.

%  The  main objective of the present paper  is to understand  to which extend the motion by mean curvature in  the classical Allen-Cahn problem holds for the system \eqref{eq1}-\eqref{eq2}  with gradient coupling and  the  volume constraint described by the Lagrange multiplier   \eqref{lagrange}. 
%  
 The paper is organized as follows.  In section \ref{well-posedness}  is devoted to the objective (i).  Section  \ref{limit} the objectives (ii) and (iii) are  addressed in the context of a model one-dimensional problem.   In Section \ref{formalderivation} the equation for the interface motion \eqref{motion1} is formally derived.

 %%%%%%%%%%%%%%%%%%%%%%%%%%%%%%%%%%%%%%%%%%%%%%%%%%%%%%%
%%%%%%%%%%%%%%%Section: Well-Posedness %%%%%%%%%%%%%%%%%%%%%%%%%%%%%
%%%%%%%%%%%%%%%%%%%%%%%%%%%%%%%%%%%%%%%%%%%%%%%%%%%%%%%

\section{Existence  of  the sharp interface solutions that do not blow up in finite time}
\label{well-posedness}

In this section we consider  the boundary value problem  \eqref{eq1}, \eqref{eq2} with $\lambda_\ve$ given by \eqref{lagrange}.  Introduce the following %\textcolor{red}{energy} 
functionals

\begin{equation}
\label{energy}
E_\ve(t):=\frac{\ve}{2} \int_\Omega |\nabla \rho_\ve(x,t)|^2dx+\frac{1}{\ve}\int_\Omega W( \rho_\ve(x,t) )dx,
\;  F_{\ve}(t):= \int_\Omega \Bigl(| P_\ve(x,t)|^2+|P_\ve(x,t)|^4\Bigr)dx.
\end{equation}
 
 The system \eqref{eq1}-\eqref{eq2} is supplied with ``well prepared" initial data, which means that two conditions hold:
 \begin{equation}
\label{rho}
-\ve^{1/4}\leq \rho_\ve (x,0)\leq 1+\ve^{1/4},
\end{equation}
 and 
 \begin{equation}
E_\ve(0)+F_\ve(0)\leq C ,
\label{IniEnBound}
\end{equation}

The first  condition \eqref{rho} is a weakened form of a standard  condition for the phase field variable $0\leq \rho_\ve(x,0)\leq 1$. If $\lambda_\ve \equiv 0$, then  by the maximum principle  $0\leq \rho_\ve(x,0)\leq 1$ implies $0\leq \rho_\ve(x,t)\leq 1$ for $t>0$.  The presence of  nontrivial  $\lambda_\ve$ leads to  an  ``extended interval"  for $ \rho_\ve$. Second condition \eqref{IniEnBound} means  
that at $t=0$ the function 
$\rho_\ve$ has the structure of ``$\ve$-transition layer", that is the domain $\Omega$  consists of three subdomains:  a subdomain 
where  the phase field function $\rho_\ve$ is close to $1$ (inside  the cell)  an
another subdomain where $\rho_\ve \sim 0$ (outside the cell) separated by a  transition region of width $\ve$. Furthermore, the  orientation field  $P_\ve$ has value close to $0$ everywhere except  the $\ve$-transition region.

\begin{theorem} If the initial data $\rho^{i}_\ve:=\rho_\ve(x,0)$, $P^i_\ve:=P_\ve(x,0)$ satisfy  \eqref{rho} and
\eqref{IniEnBound}, then for any $T>0$ the solution $\rho_\ve$, $P_\ve$ exists on the time interval 
$(0,T)$ for sufficiently small $\ve>0$, $\ve<\ve_0(T)$. 
Moreover, it satisfies  $-\ve^{1/4}\leq \rho_\ve (x,t)\leq 1+\ve^{1/4}$ and
\begin{equation}\label{noblowup1}
\ve \int_0^T\int_\Omega \Bigl(\frac{\partial\rho_\ve}{\partial t} \Bigr)^2dxdt\leq C, 
\quad 
E_\ve(t)+F_\ve(t)\leq C\quad \forall t\in(0,T),
\end{equation}
where $C$ is independent of $t$ and $\ve$. 
\end{theorem}

For the proof of this theorem we refer to \cite{BerPotRyb2016}. 

%%%%%%%%%%%%%%%%%%%%%%%%%%%%%%%%%%%%%%%%%%%%%%%%%%%%%%%
%%%%%%%%%%%%%%%   Section Sharp Interface Limit  %%%%%%%%%%%%%%%%%%%%%%%%%%
%%%%%%%%%%%%%%%%%%%%%%%%%%%%%%%%%%%%%%%%%%%%%%%%%%%%%%%
\section{1D model problem:  rigorous derivation  of the   sharp interface limit  and remarks on stability}\label{limit}

In Section \ref{formalderivation}  we present  derivation of the formal  asymptotic  expansion of the solution of \eqref{eq1}-\eqref{eq2} and use it to obtain the equation of motion \eqref{motion1}, which is the principal object of interest in the study of cell motility.  There are two main  sources of difficulties in rigorous justification of  
this derivation (i) the possible non-smoothness of  limiting velocity field  $V$ and  (ii) dimension greater than one is much harder to handle technically.
That is why in this Section we consider a simplified one-dimensional model  and  impose smallness assumption on $\beta$ that   guarantees  regularity of  asymptotic solutions  as well as plays an important technical role in our proof. {\it Recently, we have developed another approach of justification of Sharp Interface Limit for the system \eqref{orig_1}-\eqref{orig_2} which is valid for arbitrary $\beta$, see \cite{BerPotRyb2016}.}

 Specifically, we study the limiting behavior as $\varepsilon\to 0$ of the solution of the system
\begin{empheq}[left=\empheqlbrace]{align}
\frac{\partial \rho_{\varepsilon}}{\partial t}&=\partial^2_{x}\rho_{\varepsilon}-\frac{W'(\rho_{\varepsilon})}{\varepsilon^2}+P_{\varepsilon}\partial_x\rho_{\varepsilon}+\frac{F(t)}{\varepsilon},\quad x\in\mathbb{R}^1 \label{orig_1}\\ \nonumber\\ 
\frac{\partial P_{\varepsilon}}{\partial t}&=\varepsilon \partial_{x}^2P_{\varepsilon}-\frac{1}{\varepsilon}P_{\varepsilon}+\beta \partial_{x}\rho_{\varepsilon},\label{orig_2}
\end{empheq}
assuming that the initial data $\rho_{\varepsilon}(x,0)$, $P_\ve(x,0)$ has the following (``very well-prepared") form
\begin{equation}
\rho_{\varepsilon}(x,0)=\theta_0\left(\frac{x}{\ve}\right) +\sum\limits_{i=1}^{N}\ve^{i}\theta_i\left(\frac{x}{\ve}\right)+\varepsilon^{\alpha} u_{\ve}\left(\frac{x}{\ve},0\right) 
\nonumber \label{iniOneD}
\end{equation}
and 
\begin{equation}\nonumber
P_{\varepsilon}(x,0)=\sum\limits_{i=1}^{N}\ve^{i}\Psi_i\left(\frac{x}{\ve}\right)+\varepsilon^{\alpha} Q_{\ve}\left(\frac{x}{\ve},0\right),
\end{equation}
where $\alpha<N+1$. Here functions $\theta_i$ are solutions of  \eqref{theta_0} for $i=0$ and \eqref{theta_i} for $i \geq 1$, and $\Psi_i$ are solutions of  \eqref{psi0} for $i=0$ and \eqref{Psi_i} for $i\geq 1$. The functions $V_i$ which are involved in the definition of $\Psi_i$ are defined by \eqref{main_solvability} and \eqref{solvability} ($V=V_0+\ve V_1+...$  is the expansion the velocity of the cell's sharp interface). 
We also assume that there exists a constant $C$, independent of $\ve$, such that 
\begin{equation}
\|u_{\ve}(y,0)\|_{L^2}+\|Q_{\ve}(y,0)\|_{L^2}\leq C.\label{iniOneD}
\end{equation}

We emphasize that  $F(t)$ in the RHS of  \eqref{orig_1} is a given function rather than an unknown  Lagrange   multiplier in \eqref{eq1}. The main distinction of 1D case is  because in 1D there is no motion due to curvature of the interface since the  interface is a point.  Therefore  if we take initial data such that the domain $\rho =1$ is a finite interval, then  such a ``one dimensional mathematical cell" will not move, which  corresponds to a well known fact  that   the motion in the  one-dimensional   Allen-Cahn problem is exponentially slow (that is very different from AC in higher dimensions). Thus, we choose initial data to be a step like function that is a transition from   an unbounded left  interval $\rho =0$ to   an unbounded right  interval $\rho =1$ (the ``cell").  In  the 2D problem Lagrange multiplier $\lambda_\ve(t)$ appeared due to the mass (volume) conservation constraint in a finite domain occupied by the cell, which has no analog in 1D  because the ``mathematical cell" must be infinite as explained above.  Therefore an analog of the Lagrange multiplier in 1D is chosen to be a  given forcing term $F(t)$.

Hereafter  $\theta_0$ denotes the classical standing wave solution of 
\begin{equation}
\theta_0^{\prime\prime}(x)=W^{\prime}\bigl(\theta_0(x)\bigr) \quad x\in \mathbb{R}^1.
\label{deftheta_0}
\end{equation} 
with step-like conditions at infinity
\begin{equation}
\theta_0(x)\to 0,\ \text{as}\ x\to -\infty;\quad  \theta_0\to 1,\ \text{as} \ x \to +\infty.
\label{1deftheta_0}
\end{equation} 
Since the solution of \eqref{deftheta_0}-\eqref{1deftheta_0} is uniquely defined up to translations  (shifts),
we impose the following normalization
\begin{equation}
\theta_0(0)=1/2.
\label{2deftheta_0}
\end{equation} 
In the particular case of the double-well potential $W$ having the form \eqref{doublewell} the solution $\theta_0$ is explicitly given by $\theta_0 =(1-\tanh(x/\sqrt{8}))/2$. Note that $\rho_\ve(x,t)=\theta_0(x/\ve)$ solves \eqref{orig_1}  if the  coupling term and $F(t)$ are both identically zero.

The main result of this sections is the following theorem:\\
\begin{theorem} \label{theorem_sharp} Let $\rho_{\varepsilon}$ and $P_{\varepsilon}$ solve \eqref{orig_1} and \eqref{orig_2} on $[0,T]$ with initial conditions satisfying \eqref{iniOneD}. Assume also that  
$F(t)$ is a given smooth function and $|\beta| < \beta_0$, where  $\beta_0 >0$ depends on the potential $W$ only. Then we have, for sufficiently small $\ve$,
\begin{equation}\label{expansion_rho}
\rho_{\varepsilon}(x,t)=\theta_0\left(\frac{x-x_{\varepsilon}(t)}{\varepsilon}\right)+\varepsilon \rho^{(1)}_{\varepsilon}(t,x),
\end{equation}
where $x_{\varepsilon}(t)$  denotes  the location of the interface between $0$ and $1$ phases and 
  %$\varepsilon$ 
\begin{equation}\label{residual_bound}
\int (\rho^{(1)}_\ve(t,x))^2dx \leq C\ \text{ for all }\ t\in[0,T].
\end{equation}
Moreover,   $x_{\varepsilon}(t)$  converges to $x_0(t)$ which solves the interface motion 
equation similar to \eqref{motion1}: 
\begin{equation}
\label{motion1d}
-c_0\dot{x}_0(t)=\beta \Phi(\dot{x}_0(t))+F(t), \;\;\;\Phi(V):=\frac{1}{\beta}\int \Psi_0(y,-V)(\theta_0')^2dy,
\end{equation}
where $c_0=\int (\theta_0')^2 dy$, and $\Psi_0(y,V)$ is defined as solution of \eqref{def_of_psi}.
\end{theorem}

\begin{remark}  From the definition of function $\Psi_0$ it follows that it depends linearly on the parameter $\beta$ so that function $\Phi(V)$ does not depend on $\beta$.
\end {remark} 

\begin{remark}  While Theorem \ref{theorem_sharp} describes  the  leading term of the asymptotic expansion for $\rho_\ve$, in   the course of the proof we also construct the  leading term of the asymptotic expansion of $P_\ve$ in the form $\Psi_{0}\left(\frac{x-x_{0}(t)}{\varepsilon},\dot x_0(t)\right)$.
\end{remark}

\begin{remark}(on stability)  Equation \eqref{motion1d} is rigorously derived  when $|\beta| < \beta_0$ but it could be formally derived for any real  $\beta$. This equation has the unique smooth solution  ${x}_0(t)$ when $|\beta| < \beta^*$,  for some $\beta^* >\beta_0 >0$.  Roughly speaking,  if $\beta < \beta^*$,  then   ${x}_0(t)$ is determined   by  \eqref{motion1d} due to the implicit function theorem otherwise  multiple solutions ${x}_0(t)$  may appear. Thus, assumption on smallness of $\beta$ can be viewed as a stability condition. By contrast, for large enough $\beta$ one can observe instability due to the fact that the limiting equation \eqref{motion1d} has multiple solutions and, therefore, perturbation of initial data may result in switching between multiple solutions of equation \eqref{motion1d}. Indeed, to explain we rewrite equation \eqref{motion1d}
\begin{equation}\label{rewritten}
c_0V-\beta \Phi(V) =F_0(t)    
\end{equation}
where the left hand side of \eqref{rewritten} can be resolved in $V$, but not uniquely. 
\end{remark}

\begin{remark}  The estimate \eqref{residual_bound}  justifies  the asymptotic expansion \eqref{expansion_rho} and this estimate is the  principal claim of this Theorem. However, in the course of the proof we actually obtain and justify a more precise asymptotic expansion of the from
$\rho_\ve= \theta_0(\frac{x-x_\ve}{t})+\ve \theta_1(\frac{x-x_\ve}{t},t)+
\ve^2 \theta_2(\frac{x-x_\ve}{t},t)+O(\ve^3)$, which corresponds to $N=3$ and $\alpha=3$ in the the expansions \eqref{expansions_3} below.
\end{remark}

The proof of Theorem \ref{theorem_sharp}   is divided into two steps, presented in the following two subsections.

In the first step (Subsection \ref{formal_expansion}) we  formally construct approximate solution of the order $N$ and 
obtain equations for the residuals  
\begin{equation}\label{residuals}
u_{\varepsilon}=\frac{1}{\ve^\alpha}(\rho_\ve -\tilde\rho_\ve) \  \  \text{and}  \ \   
Q_{\varepsilon}=\frac{1}{\ve^\alpha} 
(P_\ve-\tilde P_\ve),
\end{equation}
 where 
\begin{equation}\label{expansions_3}
\tilde{\rho}_{\varepsilon}=\theta_0+\varepsilon\theta_1+...+\varepsilon^{N}\theta_N\text{ and }\tilde{P}_{\varepsilon}=\Psi_0+\varepsilon\Psi_1+...+\varepsilon^N\Psi_N,
\end{equation}
for some $\alpha\geq 1$. It would be natural to expect that $\alpha=N+1$, however, it turned out that due to the gradient coupling  and  nonlinearity of the  problem, we can only prove boundedness of $u_{\varepsilon}$ and $Q_{\varepsilon}$  for some  $1<\alpha < N+1$.

The second step (Subsection \ref{a_priori_bounds})  is the central mathematical part of this paper, and we briefly outline its main ideas.
The goal  there  is to  obtain bounds on residuals $u_{\varepsilon}$ and $Q_{\varepsilon}$ for appropriate $\alpha$ and $N$.   The bounds on $u_{\varepsilon}$  play central role and they imply bounds on $Q_{\varepsilon}$ though these bounds are  coupled. Therefore the bound \eqref{residual_bound} is the main  claim of the Theorem.
 
   The  techniques  of  asymptotic expansions that include  bounds on residuals  were  first developed  for  Allen-Cahn PDE in \textcolor{black}{\cite{MotSha95}}.   The proofs in  \textcolor{black}{\cite{MotSha95}} are based on the  lower bound of the  spectrum  of linearized self-adjoint stationary Allen-Cahn operator in  an unbounded domain. The techniques of this type were subsequently developed  and applied in \textcolor{black}{Alikakos, Bates, Chen [1] for the Cahn�Hillard equation, Caginalp and Chen [6] for the phase field system},  and  \cite{CheHilLog10} for  volume preserving Allen-Cahn PDE. 
   
   In the system  \eqref{eq1}-\eqref{eq2} or its one-dimensional analog    \eqref{orig_1}-\eqref{orig_2}  the  corresponding linearized operator is not self-adjoint and the previously developed techniques can not be directly applied.
   
   The results of  this Section are based on the analysis of  a time-dependent linearized problem that corresponds to the entire system. We represent (split) the  residual $u_\ve$  as a some of  two  parts  $u_{\varepsilon}(x,t)=\theta_0^\prime (x/\ve)v_\ve(x/\ve,t)+\theta_0^\prime(x/\ve) \xi_\ve(t)$, where  $v_\ve(x/\ve,t)$ and $\xi_\ve(t)$  are new unknown functions. The first part is easer to estimate since it is  chosen to be orthogonal to  the  eigenfunction  $\theta_0^\prime(x/\ve)$ corresponding  to the zero eigenvalue of the linearized  stationary Allen-Cahn operator,  but is has a more general form than the second one since $v_\ve(x/\ve,t)$ depends on both spatial and time variables.  The second part is of a simpler form  because $\xi_\ve(t)$ does not depend on $x$, but it contains the eigenfunction $\theta_0^\prime(x/\ve)$. %that corresponds to the zero eigenvalue. 
The difficulty of dealing with such an  eigenfunction can be explained by analyzing the equations \eqref{eq_for_u} and \eqref{eq_for_Q} obtained by rescaling of  the spatial variable.
Note that the two highest order $\ve^{-2}$ terms  in \eqref{eq_for_u} dominate other %lower order 
terms, and the sum of these two terms is nothing but  linearized stationary Allen-Cahn  operator.  
%whose eigenfunction $\theta_0^\prime$ corresponds to the zero (lowest) eigenvalue.  
If in \eqref{eq_for_u} one takes $u_{\varepsilon}(x,t)=\theta_0^\prime(x/\ve) \xi_\ve(t)$, then the $\ve^{-2}$ terms vanish ($\theta_0^\prime$ is an eigenfunction)  and one has  to  estimate $\xi_\ve(t)$ by analyzing   the lower  order terms, which is a much harder task. For example, in order to estimate $\xi_\ve(t)$ we   represent $Q_\ve$  in \eqref{partition_Q} as a sum of two parts  corresponding to  the representation of $u_{\varepsilon}$ and write down the leading term for the second part (first term in \eqref{repr_for_Q}). The justification of  expansion with  this leading term is a subtle task because it requires bounds on both  $\xi_\ve(t)$ and  $\dot \xi_\ve(t)$, which leads to a condition on smallness of $\beta$.

%%%%%%%%%%%%%%%%%%%%%%%%%%%%%%%%%%%%%%%%%%%%%%%%%%%%%%%
%%%%%%%%%%%%%%%%%% Formal Expansions %%%%%%%%%%%%%%%%%%%%%%%%%%%%
%%%%%%%%%%%%%%%%%%%%%%%%%%%%%%%%%%%%%%%%%%%%%%%%%%%%%%%

\subsection{Construction of  asymptotic expansions.}
\label{formal_expansion}

First, we seek  formal approximations for $\rho_{\varepsilon}$ and $P_{\varepsilon}$ in the form:
\begin{equation} \nonumber%\label{expansions_rho_P}
{\rho}_{\varepsilon}(x,t)\approx\theta_0\left(\frac{x-x_{\varepsilon}(t)}{\varepsilon}\right)+\sum\limits_{i}^{}\varepsilon^{i}\theta_i\left(\frac{x-x_\varepsilon(t)}{\varepsilon},t\right)\text{ and }%\tilde
{P}_{\varepsilon}(x,t)\approx \sum\limits_{i}^{}\varepsilon^{i}\Psi_{i}\left(\frac{x-x_{\varepsilon}(t)}{\varepsilon},t \right).
\end{equation}
We  also  assume  $x_{\varepsilon}(t)$  admits a power series expansion, 
%and $F_{\varepsilon}(t)$ are given functions and 
\begin{equation}\label{expansion_for_x}
x_{\varepsilon}(t)=x_0(t)+\varepsilon x_1(t)+...+\varepsilon^{N}x_N(t)+...,
%F_{\varepsilon}(t)&=&F_0(t)+\varepsilon F_1(t)+...+\varepsilon^{N}F_N(t)+...\nonumber
\end{equation} 
so that we also have expansion for the velocity  $V=-\dot x_\ve$,
$$
V(t)=V_0(t) +\varepsilon V_1(t)+...+\varepsilon^{N}V_N(t)+...,\quad V_i(t)=-\dot{x}_i(t), \ i=0,1,...
$$ 
 Next we expand $W'(\rho_{\varepsilon})$, 
\begin{eqnarray*}
W'(\rho_{\varepsilon})&=
%W'(\theta_0)+W''(\theta_0)\left[\varepsilon \theta_1+\varepsilon^2 \theta_2+...+\varepsilon^N\theta_N\right]\\
%&&+\frac{W'''(\theta_0)}{2}\left[\varepsilon \theta_1+\varepsilon^2 \theta_2+...+\varepsilon^N\theta_N\right]^2\\
%&&+\frac{W^{\text{(iv)}}(\theta_0)}{6}\left[\varepsilon \theta_1+\varepsilon^2 \theta_2+...+\varepsilon^N\theta_N\right]^3\\
W'(\theta_0)+\varepsilon W''(\theta_0)\theta_1+\varepsilon^2\left[W''(\theta_0)\theta_2+\frac{W'''(\theta_0)}{2}\theta_1^2\right]+...
\\ &  +\varepsilon^{i}\left[W''(\theta_0)\theta_i+(dW)^{(i)}\right]+... ,
\end{eqnarray*}
where 
\begin{equation}
(dW)^{(i)}=\sum\limits_{\footnotesize\begin{array}{c}i_1+i_2=i,\\i_1,i_2\geq 1\end{array}}\frac{W'''(\theta_0)}{2}\theta_{i_1}\theta_{i_2}+\sum\limits_{\footnotesize\begin{array}{c}i_1+i_2+i_3=i,\\i_1,i_2,i_3\geq 1\end{array}}\frac{W^{(\text{iv})}(\theta_0)}{6}\theta_{i_1}\theta_{i_2}\theta_{i_3}.\nonumber
\end{equation}
Substitute  the   series \eqref{expansions_rho_P}  and  \eqref{expansion_for_x} into \eqref{orig_1}-\eqref{orig_2}, and equate terms of like powers of $\ve$ to obtain that    
$\theta_i$  and $\Psi_i$ for $i=0,1,2$ satisfy
\begin{eqnarray}
 \theta''_0&=&W'(\theta_0)\label{theta_0}\\
-\theta_1''+W''(\theta_0)\theta_1&=&-V_0\theta_0'+\Psi_0\theta'_0+F(t)\label{theta_1}\\
-\theta_2''+W''(\theta_0)\theta_2&=&-V_1\theta'_0-V_0\theta'_1-\frac{W'''(\theta_0)}{2}\theta_1^2+\Psi_0\theta_1'+\Psi_1\theta_0'%+F_1(t)
\label{theta_2}
\end{eqnarray}
%\textcolor{red}{Add solvability conditions and properties of the operator $-\partial^2_{y}+W''(\theta_0)$}\\
and 
\begin{eqnarray}
\Psi_0''-V_0\Psi'_0-\Psi_0&=&-\beta\theta'_0,\label{psi0}\\
\Psi_1''-V_0\Psi'_1-\Psi_1&=&-\beta\theta'_1+V_1\Psi'_0+\dot{\Psi}_0,\nonumber\\
\Psi_2''-V_0\Psi'_2-\Psi_2&=&-\beta\theta'_2+V_1\Psi'_1+V_2\Psi'_0+\dot{\Psi}_1.\nonumber
\end{eqnarray}

The equations for $i>2$  have the following form
\begin{equation}
-\theta_i''+W''(\theta_0)\theta_i=-\dot\theta_{i-2}-\sum_{j=0}^{i-1}V_j\theta_{i-1-j}'-(dW)^{(i)}+\sum\limits_{j=0}^{i-1}\Psi_j\theta_{i-1-j}'.%+F_{i-1}.
\label{theta_i}
\end{equation}

\begin{equation}
\Psi_i''-V_0\Psi'_i-\Psi_i=-\beta\theta'_i+\sum\limits_{j=1}^{i}V_j\Psi'_{i-j}+\dot{\Psi}_{i-1}.\label{Psi_i}
\end{equation}

\begin{remark} \label{SOLVABILITY}
Due to the fact that $\theta^\prime$ is an eigenfunction of the linearized Allen-Cahn operator corresponding to the zero eigenvalue, the following 
solvability conditions for \eqref{theta_1},\eqref{theta_2}  and \eqref{theta_i} arise 
\begin{equation}
\int \left\{-V_0\theta_0'+\Psi_0\theta_0'+F(t)\right\}\theta_0'dy=0,\label{main_solvability}
\end{equation}
\begin{equation}\label{between_solvability}
\int \left\{-V_1\theta'_0-V_0\theta'_1-\frac{W'''(\theta_0)}{2}\theta_1^2+\Psi_0\theta_1'+\Psi_1\theta_0'\right\}\theta_0'dy=0,
\end{equation}
\begin{equation}\label{solvability}
\int \left\{-\dot\theta_{i-2}-\sum_{j=0}^{i-1}V_j\theta_{i-1-j}'-(dW)^{(i)}
+\sum\limits_{j=0}^{i-1}\Psi_j\theta_{i-1-j}'
%+F_{i-1}
\right\}\theta_0'dy =0.
\end{equation}
%Below $\theta_i$ and $\Psi_i$ are functions which are smooth, bounded and given. 
and uniquely define the functions $V_i(t)$, $i=0,...,N-1$ such that the solvability conditions \eqref{solvability}
are satisfied. Equations \eqref{main_solvability}, \eqref{between_solvability} and \eqref{solvability} are solvable for $V_0$, $V_1$ and $V_i$, respectively, for sufficiently small $\beta$. 
Also we note that once the solvability conditions are satisfied then equations \eqref{theta_1},\eqref{theta_2} and \eqref{theta_i} have a family of solutions: $\theta_i=\overline\theta_i+\gamma \theta_0'$, $\gamma\in \mathbb R$, where $\overline{\theta}_i$ is a particular solution. We choose $\gamma$ s.t. 
\begin{equation}\nonumber
\int \theta_0' \theta_i dy =0.
\end{equation}
\end{remark}

%%%%%%%%%%%%%%%%%%%%%%%%%%%%%%%%%%%%%%%%%%%%%%%%%%%%%
%%%%%%%%%%%%%%%%%Subsection: Equations for u and Q%%%%%%%%%%%%%%%%%%%%%
%%%%%%%%%%%%%%%%%%%%%%%%%%%%%%%%%%%%%%%%%%%%%%%%%%%%%
 
%\subsection{Equations for residuals $u_{\varepsilon}$ and $Q_{\varepsilon}$.}
%\label{residuals}
Define functions $u_{\varepsilon}(y,t)$ and $Q_{\varepsilon}(y,t)$ by 
\begin{equation}\label{basic_repr}
\rho_{\varepsilon}=\tilde{\rho}_{\varepsilon}(y,t)+\varepsilon^{\alpha}u_{\varepsilon}(y,t),\text{ and } \ 
P_{\varepsilon}=\tilde{P}_{\varepsilon}(y,t)+\varepsilon^{\alpha}Q_{\varepsilon}(y,t) \text{  for  }y=\frac{x-x_{\varepsilon}(t)}{\varepsilon},\end{equation}
where 
\begin{equation} \label{expansions_rho_P}
\tilde{\rho}_{\varepsilon}(y,t)=\theta_0(y)+\sum\limits_{i=1}^{N}\varepsilon^{i}\theta_i(y,t)\text{  and  } \tilde{P}_{\varepsilon}(y,t)=\sum\limits_{i=0}^{N}\varepsilon^{i}\Psi_{i}(y,t).
\end{equation}
Substituting the representation for $\rho_{\varepsilon}$ from \eqref{basic_repr}  into \eqref{orig_1} we derive the PDE for $u_{\varepsilon}$ (note that the differentiation in time and new spatial variable $y=\frac{x-x_{\varepsilon}(t)}{\varepsilon}$ are no longer independent) 
\begin{equation}\label{eq_for_u}
\frac{\partial u_{\varepsilon}}{\partial t}=\frac{u_{\varepsilon}''}{\varepsilon^2}-\frac{V_0u_{\varepsilon}'}{\varepsilon}-\frac{W''(\theta_0)u_{\varepsilon}}{\varepsilon^2}-\frac{W'''(\theta_0)u_{\varepsilon}\theta_1}{\varepsilon}+\frac{\Psi_0u_{\varepsilon}'}{\varepsilon}+\frac{Q_{\varepsilon}\theta_0'}{\varepsilon}+R_{\varepsilon}(t,y),
\end{equation} 
where $R_{\varepsilon}$ is of the form
\begin{eqnarray}
R_{\varepsilon}(t,y)&=&\varepsilon^{N-1-\alpha}a_{\varepsilon}(t,y)\nonumber\\&&+\varepsilon^{N-\alpha}b_{0,\varepsilon}(t,y)+b_{1,\varepsilon}(t,y)u_{\varepsilon}
+\varepsilon^{\alpha-2}b_{2,\varepsilon}(t,y)u_{\varepsilon}^2+\varepsilon^{2\alpha-2}b_{3,\varepsilon}(t,y)u_{\varepsilon}^3\nonumber\\
&&+e_{\varepsilon}(t,y)u_{\varepsilon}'+g_{\varepsilon}(t,y)Q_{\varepsilon}+\varepsilon^{\alpha-1}Q_{\varepsilon}u_{\varepsilon}'.\label{repr_for_R_eps}
\end{eqnarray}
where $a_\varepsilon(t,y),$ $b_{k,\varepsilon}(t,y),k=1,2,3$, $e_{\varepsilon}(t,y),g_{\varepsilon}(t,y)$ are bounded functions in $y$, $t$ and $\varepsilon$ and square integrable  with respect to $y$ (except $e_{\varepsilon}$). Moreover, the function $a_\ve$ is orthogonal to $\theta_0'$: 
\begin{equation}\nonumber
\int \theta_0'(y) a_\ve(t,y)dy=0. 
\end{equation}
The functions $a_\ve$, $b_{i,\ve}$, $e_{\ve}$ and $g_{\ve}$ are expressed in terms of $\theta_i$, $V_i$ and $\Psi_i$, their exact 
form, which is not important for the proof,  is given in Appendix. 

Substituting the representation  for $P_{\varepsilon}$ from \eqref{basic_repr} into \eqref{orig_2} we derive also the 
PDE for 
$Q_{\varepsilon}$:
\begin{equation}\label{eq_for_Q}
\frac{\partial Q_{\varepsilon}}{\partial t}=\frac{Q_{\varepsilon}''}{\varepsilon}-\frac{VQ_{\varepsilon}'}{\varepsilon}-\frac{Q_{\varepsilon}}{\varepsilon}+ \frac{\beta u_{\varepsilon}'}{\varepsilon}+\varepsilon^{N-\alpha}m_{\varepsilon}(t,y).
\end{equation}
%All introduced functions - $a_\varepsilon(t,y),$ $b_{k,\varepsilon}(t,y),k=1,2,3$, $e_{\varepsilon}(t,y),g_{\varepsilon}(t,y)$, $m_{\varepsilon}(t,y)$ - are given in Appendix.
For more details on derivation of \eqref{eq_for_u} and \eqref{eq_for_Q} we refer to Appendix.

%%%%%%%%%%%%%%%%%%%%%%%%%%%%%%%%%%%%%%%%%%%%%%%%%%%%%%%
%%%%%%%%%%%%%%%%%%Section: A priori bounds %%%%%%%%%%%%%%%%%%%%%%%%%%
%%%%%%%%%%%%%%%%%%%%%%%%%%%%%%%%%%%%%%%%%%%%%%%%%%%%%%%

\subsection {Bounds for residuals $u_\ve$ and $Q_\ve$}
\label{a_priori_bounds}
In this section we obtain bounds for the coupled system of PDEs: 
\begin{empheq}[left=\empheqlbrace]{align}\frac{\partial u_{\varepsilon}}{\partial t}&=\frac{u_{\varepsilon}''}{\varepsilon^2}-\frac{V_0u_{\varepsilon}'}{\varepsilon}-\frac{W''(\theta_0)u_{\varepsilon}}{\varepsilon^2}-\frac{W'''(\theta_0)u_{\varepsilon}\theta_1}{\varepsilon}+\frac{\Psi_0u_{\varepsilon}'}{\varepsilon}+\frac{Q_{\varepsilon}\theta_0'}{\varepsilon}+R_{\varepsilon}(t,y)\nonumber\\ 
\frac{\partial Q_{\varepsilon}}{\partial t}&=\frac{Q_{\varepsilon}''}{\varepsilon}-\frac{VQ_{\varepsilon}'}{\varepsilon}-\frac{Q_{\varepsilon}}{\varepsilon}+ \frac{\beta u_{\varepsilon}'}{\varepsilon}+\varepsilon^{N-\alpha}m_{\varepsilon}(t,y).\nonumber
\end{empheq}
To this end we write the unknown function $u_{\varepsilon}$ in the following form,
\begin{equation}
u_{\varepsilon}(t,y)=\theta'_0(y)\left[v_{\varepsilon}(t,y)+\xi_{\ve}(t)\right],\text{ where }\int (\theta_0'(y))^2v_{\varepsilon}(t,y)dy=0.\label{def_of_v}
\end{equation}
Then \eqref{eq_for_u}  becomes
\begin{eqnarray}\label{eq_for_fact_u}
\noindent\frac{\partial}{\partial t}\left(\theta'_0(v_{\varepsilon}+\xi_{\ve})\right)&=&-\frac{V_0}{\varepsilon}(\theta_0'(v_{\varepsilon}+\xi_{\ve}))'+\frac{(\theta_0'(v_{\varepsilon}+\xi_{\ve}))''}{\varepsilon^2}\nonumber\\\nonumber
\noindent&&-\frac{W''(\theta_0)}{\varepsilon^2}\theta_0'(v_{\varepsilon}+\xi_{\ve})-\frac{W'''(\theta_0)}{\varepsilon}(v_{\varepsilon}+\xi_{\ve})\theta_0'\theta_1\\&&+\frac{\Psi_0(\theta_0'(v_{\varepsilon}+\xi_{\ve}))'}{\varepsilon}+\frac{Q_{\varepsilon}\theta_0'}{\varepsilon}+R_{\varepsilon}(t,y). 
\end{eqnarray}

\begin{lemma}  The following inequality holds
\begin{eqnarray}
&&\frac{d}{2dt}\left[\int (\theta_0')^2v_{\varepsilon}^2dy+c_0 \xi_{\ve}^2\right]+\frac{1}{2\varepsilon^2}\int (\theta_0')^2(v'_{\varepsilon})^2dy \nonumber\\
&&\hspace{-10 pt}\leq C \xi_{\ve}^2+\frac{1}{\varepsilon}\left[\int Q_{\varepsilon}(\theta_0')^2(v_{\varepsilon}+\xi_{\ve})dy-
 \xi_{\ve}^2 \int \Psi'_0(\theta_0')^2dy \right]+\int R_{\varepsilon}\theta_0'(v_{\varepsilon}+\xi_{\ve})dy.\label{en_eq_75}
\end{eqnarray}

\end{lemma}
\noindent{P r o o f:}\\
Multiply \eqref{eq_for_fact_u} by $u_{\varepsilon}=\theta_0'(v_{\varepsilon}+\xi_{\ve})$ and integrate to obtain 
\begin{eqnarray}
&&\frac{d}{2dt}\left[\int (\theta_0')^2(v_{\varepsilon}+\xi_{\ve})^2dy\right]+\frac{1}{\varepsilon^2}\int \left\{(\theta_0'(v_{\varepsilon}+\xi_{\ve}))'\right\}^2dy \nonumber\\
&&=-\frac{V_0}{\varepsilon}\int (\theta_0'(v_{\ve}+\xi_{\ve}))'\theta_0'(v_{\ve}+\xi_{\ve})dy\nonumber
\\&&\hspace{10 pt}-\int\frac{W''(\theta_0)}{\ve^2}(\theta_0')^2(v_{\ve}+\xi_{\ve})^2dy-\int \frac{W'''(\theta_0)}{\ve}(\theta_0')^2\theta_1(v_{\ve}+\xi_{\ve})^2dy\nonumber\\
&&\hspace{10 pt}+\int \frac{\Psi_0}{\varepsilon} (\theta_0'(v_{\ve}+\xi_{\ve}))'\theta_0(v_{\ve}+\xi_{\ve})dy\nonumber
\\&& \hspace{10 pt} +\int \frac{Q_{\ve}}{\varepsilon} (\theta_0')^2(v_{\ve}+\xi_{\ve})dy+\int R_{\varepsilon}\theta_0'(v_{\varepsilon}+\xi_{\ve})dy.\label{eq_75}
\end{eqnarray}

In order to derive  \eqref{en_eq_75} we simplify equality \eqref{eq_75}. 
%We will derive inequality \eqref{en_eq_75} from \eqref{eq_75}.
{First}, we notice that the integral in the first term can be rewritten as follows 
\begin{equation}\label{observ_1}
\int (\theta_0')^2(v_{\varepsilon}+\xi_{\ve})^2dy=\int (\theta_0')^2v_{\varepsilon}^2dy+\xi_{\ve}^2\int (\theta_0')^2dy.
\end{equation}
We used here the definition of $v_{\ve}$, i.e., $\int (\theta_0')^2 v_{\ve}dy=0$.

The first term in the right hand side of \eqref{eq_75} vanishes. Indeed, 
\begin{equation*}
\int (\theta_0'(v_{\ve}+\xi_{\ve}))'\theta_0'(v_{\ve}+\xi_{\ve})dy=\left.\frac{1}{2}(\theta_0')^2(v_{\ve}+\xi_{\ve})^2\right|_{y=-\infty}^{+\infty}=0.
\end{equation*}

The second term in the left hand side can be split into three terms: 
\begin{eqnarray}
\frac{1}{\varepsilon^2}\int \left\{(\theta_0'(v_{\varepsilon}+\xi_{\ve}))'\right\}^2dy&=&\frac{1}{\ve^2}\int \left((\theta_0'v_{\varepsilon})'\right)^2 dy\nonumber\\&&+\frac{2 \xi_{\ve}}{\ve^2}\int (\theta_0' v_{\ve})' \theta_0'' dy +\frac{\xi_{\ve}^2}{\ve^2}\int (\theta_0'')^2dy.\label{eq_76} 
\end{eqnarray}

%Next we combine these three  terms in the right hand side of \eqref{eq_76} with parts of the second and the third terms in the right hand side of \eqref{eq_75}. 

Rewrite the first term in \eqref{eq_76}, using \eqref{theta_0} and integrating by parts,
\begin{eqnarray}
\int\left((\theta_0'v_{\varepsilon})'\right)^2 dy&=&\int (\theta_0'')^2v_{\varepsilon}^2dy+2\int \theta_0'v_{\varepsilon}'\theta_0''v_{\varepsilon}dy+\int (\theta_0')^2(v_{\varepsilon}')^2dy\nonumber\\
&=&\left\{-\int \theta_0'\theta_0'''v_{\varepsilon}^2dy-2\int\theta_0''\theta_0'v_{\varepsilon}'v_{\varepsilon}dy\right\}+2\int \theta_0'v_{\varepsilon}'\theta_0''v_{\varepsilon}dy\nonumber\\
&&+\int (\theta_0')^2(v_{\varepsilon}')^2dy\nonumber\\
&=&-\int \theta_0'\theta_0'''v_{\varepsilon}^2dy+\int (\theta_0')^2(v_{\varepsilon}')^2dy\nonumber\\
&=&-\int W''(\theta_0)(\theta_0')^2v_{\varepsilon}^2dy+\int (\theta_0')^2(v_{\varepsilon}')^2dy.\nonumber
\end{eqnarray}
Thus,
\begin{equation}\label{observ_2}
-\frac{1}{\varepsilon^2}\int\left((\theta_0'v_{\varepsilon})'\right)^2 dy-\int \frac{W''(\theta_0)}{\varepsilon^2}(\theta_0')^2v_{\varepsilon}^2dy=-\frac{1}{\varepsilon^2}\int (\theta_0')^2(v_{\varepsilon}')^2dy.
\end{equation}

Next we rewrite the second and the third terms  in \eqref{eq_76} using the fact that 
$\theta_0'''=W''(\theta_0)\theta_0'$  (this latter equality is obtained by differentiating \eqref{theta_0} with respect to $y$),
\begin{eqnarray}
-\frac{1}{\varepsilon^2}\int (\theta_0'')^2dy\xi_{\ve}^2-\int \frac{W''(\theta_0)}{\varepsilon^2}(\theta_0'^2)dy \xi_{\ve}^2&=&0,\label{observ_4}\\
-\frac{2}{\varepsilon^2}\int (\theta_0'v_{\varepsilon})'\theta_0''dy \xi_{\ve}-2 \int \frac{W''(\theta_0)}{\varepsilon^2}(\theta_0')^2v_{\varepsilon}dy \xi_{\ve}&=&0.\label{observ_5}
\end{eqnarray}

Also, we make use the following equality which follows from \eqref{theta_2},
\begin{eqnarray}
-\int \frac{W'''(\theta_0)}{\varepsilon}\theta_1(\theta_0')^2dy+\int\frac{\Psi_0}{2\varepsilon}\left((\theta_0')^2\right)'dy &=& \int \left[-V_0 \theta_0' +\Psi_0\theta_0' +\theta_1''\right]\theta_0''dy\nonumber\\
&&+\int W''(\theta_0)\theta_0'\theta_1'dy+\int \Psi_0\theta_0'\theta_0''dy\nonumber\\&=&-\int \Psi_0'(\theta_0')^2dy.\label{observ_3}
\end{eqnarray}

Finally,  we use equalities \eqref{observ_1},\eqref{observ_2},\eqref{observ_3},\eqref{observ_4} and \eqref{observ_5} to rewrite \eqref{eq_75} as follows
\begin{eqnarray}
&\frac{d}{2dt}\left[\int (\theta_0')^2v_{\varepsilon}^2dy+\int (\theta_0')^2dy \xi_{\ve}^2\right]&+\frac{1}{\varepsilon^2}\int (\theta_0')^2(v'_{\varepsilon})^2dy=\nonumber\\
&&-2\xi_{\ve} \int \frac{W'''(\theta_0)}{\varepsilon}\theta_1(\theta_0')^2v_{\varepsilon}dy-\int\frac{W'''(\theta_0)}{\varepsilon}\theta_1(\theta_0')^2v_{\varepsilon}^2dy\nonumber\\&&-\xi_{\ve} \int \frac{\Psi_0'}{\varepsilon}(\theta_0')^2v_{\varepsilon}dy-\int \frac{\Psi'_0}{2\varepsilon}(\theta_0')^2v_{\varepsilon}^2dy\nonumber\\
&&+\frac{1}{\varepsilon}\int Q_{\varepsilon}(\theta_0')^2(v_{\varepsilon}+\xi_{\ve})dy-\frac{\xi_{\ve}^2}{\varepsilon}\int \Psi'_0(\theta_0')^2dy \nonumber\\
&&+\int R_{\varepsilon}\theta_0'(v_{\varepsilon}+\xi_{\ve})dy.\nonumber
\end{eqnarray}
Then  \eqref{en_eq_75} is obtained by applying the standard Cauchy-Schwarz inequality and the Poincar\'e inequality \eqref{poincare} from Appendix.\\
\rightline{$\square$}

\noindent It follows from the previous lemma that in order show the boundedness of  $u_{\varepsilon}$  we need to  find an appropriate upper bound  on the term 
\begin{equation}\label{Bad_term}
\frac{1}{\varepsilon}\left[\int Q_{\varepsilon}(\theta_0')^2(v_{\varepsilon}+\xi_{\ve})dy-\xi_{\ve}^2 \int \Psi'_0(\theta_0')^2dy \right].
\end{equation}
To this end, we use  the equation \eqref{eq_for_Q} for $Q_{\varepsilon}$.

% It follows from \eqref{psi0} function $\Psi_0$ is a function of spatial variable $y$ and time $t$. Since the right hand side of \eqref{psi0}, $-\beta\theta_0'$, does not depend on time $t$, we can define $\Psi_0$ as a function of $y$ and velocity $V_0$. More precisely, let $\Psi_0(y;V_0)$  be the solution of the following equation:
%\begin{equation}\label{def_of_psi}
%\Psi_0''-V_0\Psi_0'-\Psi_0=-\beta\theta_0'.
%\end{equation}

Note that $\Psi_0=\Psi_0(y;V_0)$ depends on time $t$ through $V_0(t)$. For simplicity of the presentation we suppress the second argument $\Psi_0(y):=\Psi_0(y;V_0)$, if it equals to $V_0=-\dot{x}_0(t)$. 

It follows from \eqref{psi0} functions $\Psi_0'$ and $\Psi_{0,V_0}:=\partial \Psi_0/\partial V_0$ solve the following equations
 \begin{equation}\nonumber
(\Psi'_0)''-V_0(\Psi'_0)'-(\Psi_0')=-\beta\theta_0''\text{ and }\Psi''_{0,V_0}-V_0\Psi'_{0,V_0}-\Psi_{0,V_0}=\Psi'_0.
\end{equation}
Rewrite \eqref{eq_for_Q} substituting $u_{\ve}=\theta_0'  (v_{\ve}+\xi_{\ve})$,
\begin{equation}\nonumber
Q_{\varepsilon}''-VQ_{\varepsilon}'-Q_{\varepsilon}-\varepsilon \frac{\partial Q_{\varepsilon}}{\partial t}=-\beta \theta_0''\xi_{\ve}-\beta(\theta_0'v_{\varepsilon})'-\varepsilon^{N+1-\alpha} m_{\varepsilon}.
\end{equation}
Thus, $Q_{\ve}$ can be written as 
\begin{equation}\label{partition_Q}
Q_{\varepsilon}= \Qcc+\Qs,
\end{equation}
where $\Qcc$ and $\Qs$ are solutions of  the following problems, 
\begin{eqnarray}\label{forA}
\Qcc''-V\Qcc'-\Qcc-\varepsilon \frac{\partial\Qcc}{\partial t}&=&-\beta\theta_0''\xi_{\ve},\;\;\Qcc(0)=Q_{\varepsilon}(0)\\
\Qs''-V\Qs'-\Qs-\varepsilon \frac{\partial \Qs}{\partial t}&=&-\beta(\theta_0'v_{\varepsilon})'-\varepsilon^{N+1-\alpha} m_{\varepsilon},\;\;\Qs(0)=0.\nonumber
\end{eqnarray}
Next we note that the function $\Qcc$ can be written as $\Qcc=\xi_{\ve}\Psi'_0+\Qc$ ($\Psi_0=\Psi_0(y;V)$), where $\Qc$ solves the following problem,
\begin{equation}\label{forD}
\Qc''-V\Qc'-\Qc-\varepsilon \frac{\partial \Qc}{\partial t}=\varepsilon\dot{\xi}_{\ve}\Psi_0'+\varepsilon\dot{V}\xi_{\ve}\Psi_{0,V},\;\;\Qc(0)=Q_{\ve}(0)-\left.\xi_{\ve}\Psi'_0\right|_{t=0}.
\end{equation}
Thus,
\begin{equation}
Q_{\varepsilon}=\xi_{\ve}\Psi'_0(y;V)%\varepsilon\xi_{\ve}\left\{\frac{\Psi_0'(y;V)-\Psi'_0(y;V_0)}{\varepsilon}\right\}
+\Qc+\Qs.\label{repr_for_Q}
\end{equation} 
This representation  \eqref{repr_for_Q} allows us to rewrite the  term \eqref{Bad_term} as follows,
\begin{equation}\label{nrepr_Bad_term}
\frac{\xi_{\ve}}{\ve}\int Q_{\ve}(\theta_0')^2 v_{\ve}dy +\frac{1}{\ve}\int \left[\Qc+\Qs\right](\theta_0')^2 dy .
\end{equation}
The bounds for these terms are collected in the next lemma.

\begin{lemma} 
\label{Rem_in_Q}
The following inequalities hold
\begin{list}{}{}
\item{(i)}
\begin{eqnarray}\label{eq_for_Qs}
\nonumber&&\varepsilon\frac{d}{dt}\left[\int\Qs^2dy\right]+\int \Qs^2 dy +\int (\Qs')^2 dy\\&&\hspace{130 pt}\leq c\int(\theta_0')^2 (v_{\varepsilon}')^2 dy+c\varepsilon^{2(N+1-\alpha)};
\end{eqnarray}
\item{(ii)} 
\begin{equation}\label{en_eq_for_Qc}
\ve\frac{d}{dt}\left[\int\Qc^2dy\right]+\int \Qc^2 dy +\int (\Qc')^2 dy \leq c\beta^2\varepsilon^2 \dot{\xi}_{\ve}^2+c\varepsilon^2\xi_{\ve}^2;
\end{equation}
\item{(iii)}
\begin{equation}\label{general_Q}
\int Q_{\varepsilon}^2 dy \leq C\xi_{\ve}^2+\int \Qc^2dy + \int \Qs^2 dy;
\end{equation}
and 
\begin{equation}\label{general_Q_2}
\int (Q_{\varepsilon}')^2 dy \leq C\xi_{\ve}^2+\int (\Qc')^2 dy+\int (\Qs')^2 dy.
\end{equation}
\end{list}
\end{lemma}
\noindent{P r o o f.}\\
Items  ($i$) and ($ii$) are proved by means of energy relations that are obtained after multiplying  \eqref{forA} and \eqref{forD} by $\Qs$ and $\Qcc$, respectively, and integrating in $y$. The resulting energy relation for the function $\Qs$ is
\begin{eqnarray}
\frac{\varepsilon}{2}\frac{d}{dt}\left[\int\Qs^2dy\right]+\int \Qs^2 dy +\int (\Qs')^2 dy &=&-\beta\int \theta_0' v_{\varepsilon}\Qs'dy+\ve^{N+1-\alpha}\int m_{\varepsilon}\Qs dy,\nonumber
\end{eqnarray}
%Note that we can multiply the inequality above by $1/\varepsilon$ and, if $\beta<1$, then by ${1}/{\varepsilon^2}$.
%Thus, multiplying by $1/c_p\beta^2\varepsilon^2$ we get: 
%\begin{equation}\nonumber
%\frac{1}{2c_p\beta^2\varepsilon}\frac{d}{dt}\left[\int\Qs^2dy\right]+\frac{1}{c_p\beta^2\varepsilon^2}\left[\int \Qs^2 dy +\frac{1}{2}\int (\Qs')^2 dy\right]\leq \frac{1}{2\varepsilon^2}\int(\theta_0')^2 (v_{\varepsilon}')^2 dy+\int\frac{m_{\varepsilon}\Qs}{c_p\beta^2\varepsilon}  dy\label{eq_for_Qs}
%\end{equation}
and the energy relation for  $\Qc$ reads
\begin{eqnarray}\label{anotherenergyrelation} 
\frac{\varepsilon}{2}\frac{d}{dt}\left[\int\Qc^2dy\right]+\int \Qc^2 dy +\int (\Qc')^2 dy &=&-\varepsilon\dot{\xi}_{\ve}\int \Psi_0'\Qc dy-\varepsilon \dot{V}\xi_{\ve}\int\Psi_{0,V}\Qc dy. 
%&\leq& \varepsilon^2\frac{1}{2}\left\{\int(\Psi_0')^2dy\right\} (\dot{\xi}_{\ve})^2+C\varepsilon^2 \xi_{\ve}^2+\frac{3}{4}\int \Qc^2 dy. 
\end{eqnarray}
%Thus, 
%\begin{equation}\nonumber
%\ve\frac{d}{dt}\left[\int\Qc^2dy\right]+\int \Qc^2 dy +\int (\Qc')^2 dy \leq c\beta^2\varepsilon^2 \dot{\xi}_{\ve}^2+c\varepsilon^2\xi_{\ve}^2.
%\end{equation}
%We used that $\int(\Psi'_0)^2dy\leq c_0\beta^2$.
Then \eqref{eq_for_Qs} and \eqref{en_eq_for_Qc} are obtained by applying the Cauchy-Schwarz inequality. Note that $\Psi_0$ depends linearly on $\beta$, this allows us to bound  the right hand side of  \eqref{anotherenergyrelation} by  $c\beta^2\ve^2\dot\xi_\ve^2+
c\ve^2\xi_\ve^2+\frac{1}{2}\int \Qc^2 dy$ and this eventually leads to  \eqref{en_eq_for_Qc} .
% This can be explained by the fact that the right hand side of  linear differential equation \eqref{def_of_psi} depends linearly on $\beta$.

Item ($iii$) easily follows from representation \eqref{repr_for_Q}.\\
\rightline{$\square$}

Using \eqref{en_eq_75}, representation  \eqref{nrepr_Bad_term} for  \eqref{Bad_term} and the previous Lemma we derive the following corrolary.

\begin{corrolary}
The following inequality holds 
\begin{eqnarray}
&&\frac{d}{2dt}\left[\int (\theta_0')^2 v_{\ve}^2 dy +c_0\xi^2_{\ve}+\frac{1}{\ve}\int \Qc^2 dy +\frac{1}{\ve}\int \Qs^2 dy \right]\nonumber\\
&&\hspace{40 pt}+\frac{1}{2\ve^2}\left[\int (\theta_0')^2(v'_{\ve})^2 dy +\int \Qc^2 dy +\int (\Qc')^2 dy +\int \Qs^2 dy +\int (\Qs')^2 dy\right]\nonumber\\
&&\hspace{80 pt}\leq c\xi^2_{\ve}+ c\beta^2 \dot{\xi}^2_{\ve} +\int R_{\ve}\theta_0'(v_{\ve}+\xi_{\ve})dy+c\varepsilon^{2(N-1-\alpha)}.\label{cor_1}
\end{eqnarray}
\end{corrolary} 

Thus, we reduced the estimation of \eqref{Bad_term} to estimation of $\dot{\xi}_{\ve}^2$.

%It follows from Lemma \ref{Rem_in_Q} that if one substitute \eqref{repr_for_Q} 
% in  \eqref{Bad_term} then 
% in fact of the order 1$, however its upper bound involves both $\xi_{\ve}$ and $\dot\xi_{\ve}$.
% The function $\dot\xi_{\ve}$ is estimated in the following  

%In the following lemma we estimate $\dot \xi^2_{\ve}$ using equation \eqref{eq_for_fact_u}. 
The key observation leading to the desired bound on $\xi_{\ve}$ can be explained now as follows. Observe that the presence of 
$\xi_{\ve}$ in the RHS of \eqref{en_eq_75}  might result in the exponential growth of $\xi_{\ve}$ (the best one can guarantee from $\dot{\xi}^2_{\ve} \leq \frac{C}{\ve} \xi^2_{\ve}$ type  bound).  
Fortunately, 
the term $\int \Psi'_0(\theta_0')^2dy \xi^2_{\ve}$  in \eqref{en_eq_75} cancels  with the leading term 
appearing
% term containing  $\xi_{\ve}$ 
after substitution of  expansion  \eqref{repr_for_Q}  for $Q_\ve$. However, this results in the 
appearance of lower order terms depending on $\dot\xi_{\ve}$. 
Lemma \ref{beta_small} below provides the control of $| \dot\xi_{\ve} |$.

%Substitute \eqref{repr_for_Q} into the RHS of 
%\eqref{eq_for_fact_u}. Then
 %is that the leading term in the expansion \eqref{repr_for_Q} of $Q_\ve$ which appears in 
 %\eqref{eq_for_fact_u} with the factor $\frac{1}{\ve}$ cancel  
%other terms of order $\frac{1}{\ve}$ when one 
%multiplies \eqref{eq_for_fact_u} by $\dot{\xi}_{\ve} \theta_0^\prim$
%and integrate the resulting equation.   

\begin{lemma} 
\label{beta_small}
%If $\beta$ is sufficiently small 
The following inequality holds
\begin{eqnarray}\label{eq_for_dot_xi}
\dot{\xi}^2_{\ve}&\leq& \frac{C}{\varepsilon^2}\int (\theta_0)^2(v'_{\varepsilon})^2dy+\frac{C}{\varepsilon^2}\int \left\{\Qs^2+\Qc^2\right\}dy+C\xi^2_{\ve}+\int R_{\varepsilon}\theta_0'dy \dot{\xi}_{\ve}.
\end{eqnarray}
\end{lemma}
\noindent{P r o o f:}\\
Multiply equation \eqref{eq_for_fact_u} by $\theta_0'\dot \xi_{\ve}$ and integrate in $y$ to obtain 
\begin{eqnarray}
\nonumber \int (\theta_0')^2(\dot v_{\varepsilon}+\dot{\xi}_{\ve})\dot\xi_{\ve} dy &=& -\frac{V_0\dot{\xi}_{\ve}}{\varepsilon}\int (\theta_0'(v_{\varepsilon}+\xi_{\ve})')'\theta_0'dy \\
\nonumber&&+\frac{\dot{\xi}_{\ve}}{\varepsilon^2}\int (\theta_0'(v_{\varepsilon}+\xi_{\ve}))''\theta_0'dy \\
\nonumber&&-\dot{\xi}_{\ve}\int \frac{W''(\theta_0)}{\varepsilon^2}(v_{\varepsilon}+\xi_{\ve})(\theta_0')^2dy \\
\nonumber&&-\dot{\xi}_{\ve}\int\frac{W'''(\theta_0)}{\varepsilon}(v_{\varepsilon}+\xi_{\ve})(\theta_0')^2\theta_1dy\\&&+\dot{\xi}_{\ve}\int \frac{\Psi_0 }{\varepsilon}(\theta_0'(v_{\varepsilon}+\xi_{\ve}))'\theta_0'dy+\dot{\xi}_{\ve} \int \frac{Q_{\varepsilon}}{\varepsilon}(\theta_0')^2dy \nonumber\\
&&+\dot{\xi}_{\ve}\int R_{\varepsilon}(t,y)\theta_0' dy .\label{beta1}
\end{eqnarray}
Using \eqref{def_of_v} we simplify the left hand side of \eqref{beta1},
\begin{eqnarray}\label{simplf_beta_0}
\int(\theta_0')^2(\dot{v}_{\varepsilon}+\dot{\xi}_{\ve})\dot{\xi}_{\ve}dy=\dot{\xi}_{\ve}\frac{\partial}{\partial t}\left\{\int(\theta_0')^2v_{\varepsilon}\right\}+\dot{\xi}_{\ve}^2 \int (\theta_0')^2dy =\dot{\xi}_{\ve}^2 \int (\theta_0')^2dy .
\end{eqnarray}
In the next four steps inequality \eqref{eq_for_dot_xi} will be derived by estimating the right hand side of  \eqref{beta1} term by term.

\noindent STEP 1. The first term in the right hand side of \eqref{beta1} is estimated by using integration by parts and the Cauchy-Schwarz inequality,
\begin{eqnarray}\nonumber
&&-\frac{V_0}{\varepsilon}\int (\theta_0'(v_{\varepsilon}+\xi_{\ve}))'\theta_0'\dot{\xi}_{\ve}dy=\frac{V_0}{\varepsilon}\int \theta_0'\theta_0''(v_{\varepsilon}+\xi_{\ve})\dot{\xi}_{\ve}dy=-\frac{V_0 \dot{\xi}_{\ve}}{2\varepsilon}\int (\theta_0)^2v'_{\varepsilon}dy \\&&\hspace{100 pt}\leq \frac{C}{\delta\varepsilon^2}\int (\theta_0')^2(v'_{\varepsilon})^2dy +\delta \dot{\xi}_{\ve}^2.\label{simplf_beta_1}
\end{eqnarray}
Here we introduced small parameter $\delta$ which does not depend on $\varepsilon$ and will be chosen later. 

\noindent STEP 2. In this step we show that the sum of the second and the third terms in the right hand side of \eqref{beta1} gives zero. Indeed, use integration by parts and the definition of $\theta_0'$ ($\theta''=W'(\theta_0)$) to obtain
\begin{eqnarray}
&&\frac{1}{\varepsilon^2}\left\{\int (\theta_0'(v_{\varepsilon}+\xi_{\ve}))''\theta_0'\dot{\xi}_{\ve}dy-\int W''(\theta_0)(v_{\varepsilon}+\xi_{\ve})(\theta_0')^2\dot{\xi}_{\ve}dy\right\}\nonumber\\&&\hspace{25 pt}=\frac{\dot{\xi}_{\ve}}{\varepsilon^2}\left\{-\int(\theta_0'(v_{\varepsilon}+\xi_{\ve}))'\theta_0''+\int W'(\theta_0)(\theta_0'(v_{\varepsilon}+\xi_{\ve}))'dy\right\} =0.\label{simplf_beta_2}
\end{eqnarray}

\noindent STEP 3. In this step we estimate the sum of the fourth, the fifth and the sixth terms in the right hand side of \eqref{beta1},
\begin{equation}
-\dot{\xi}_{\ve} \int\frac{W'''(\theta_0)}{\varepsilon}(v_{\varepsilon}+\xi_{\ve})(\theta_0')^2\theta_1dy+\dot{\xi}_{\ve} \int \frac{\Psi_0 }{\varepsilon}(\theta_0'(v_{\varepsilon}+\xi_{\ve}))'\theta_0'dy+ \dot{\xi}_{\ve} \int \frac{Q_{\varepsilon}}{\varepsilon}(\theta_0')^2dy \nonumber
\end{equation}
 
To this end we first show that 
\begin{eqnarray}\label{simplification_of_beta}
\xi_{\ve}\dot{\xi}_{\ve} \int\frac{W'''(\theta_0)}{\varepsilon}(\theta_0')^2\theta_1dy+
\xi_{\ve}\dot{\xi}_{\ve} \int \frac{\Psi_0 }{\varepsilon}\theta_0''\theta_0'dy&=&0.
\end{eqnarray}
Indeed,  differentiating \eqref{theta_1} in $y$ one obtains the  equation $\theta_1'''-W''(\theta_0)\theta_0'=
W'''(\theta_0)\theta_0'\theta_1+V_0\theta_0''-\Psi'_0\theta_0'-\Psi_0\theta_0''+F'(t)$, whose solvability condition reads 
\begin{eqnarray}\nonumber
\int \left\{W'''(\theta_0)\theta_0'\theta_1+V_0\theta_0''-\Psi'_0\theta_0'-\Psi_0\theta_0''+F'(t)\right\}\theta_0'dy=0.
\end{eqnarray}
The latter equality contains five terms. The second and the fifth term vanish since they are integrals of derivatives ($(\frac{V_0}{2}(\theta_0')^2)'$ and $(F'(t)\theta_0')'$, respectively). Integrating by parts  the third term we get, 
\begin{equation}\nonumber
\int W'''(\theta_0)(\theta_0')^2\theta_1+\Psi_0\theta_0'\theta_0''dy=0.
\end{equation}
This immediately implies \eqref{simplification_of_beta}. Now, taking into account the equality \eqref{simplification_of_beta} and  representation \eqref{repr_for_Q} for $Q_{\ve}$ written in the form
\begin{equation*}
Q_{\ve}=\xi_{\ve}\Psi_0'+\xi_{\ve}(\Psi'_0(y;V)-\Psi_0(y;V_0))+\Qc+\Qs,
\end{equation*}
we get
\begin{eqnarray}
&&-\dot{\xi}_{\ve}\int\frac{W'''(\theta_0)}{\varepsilon}(v_{\varepsilon}+\xi_{\ve})(\theta_0')^2\theta_1dy+\dot{\xi}_{\ve}\int \frac{\Psi_0 }{\varepsilon}(\theta_0'(v_{\varepsilon}+\xi_{\ve}))'\theta_0'dy+ \dot{\xi}_{\ve}\int \frac{Q_{\varepsilon}}{\varepsilon}(\theta_0')^2dy\nonumber\\
&&\hspace{40 pt}=-\dot{\xi}_{\ve}\int\frac{W'''(\theta_0)}{\varepsilon}v_{\varepsilon}(\theta_0')^2\theta_1dy+\dot{\xi}_{\ve}\int \frac{\Psi_0 }{\varepsilon}(\theta_0'v_{\varepsilon})'\theta_0'dy\nonumber\\
&&\hspace{45 pt}+\dot{\xi}_{\ve}\int {(\Psi_0'(y,V)-\Psi'_0(y,V_0))}(\theta_0')^2dy+\dot\xi_{\ve}\int \frac{\Qs}{\varepsilon}(\theta_0')^2dy+\dot{\xi}_{\ve}\int \frac{\Qc}{\varepsilon}(\theta_0')^2dy,\nonumber
\end{eqnarray}
where we have also used integration by parts.

Applying the Poincar\'e inequality \eqref{poincare} we get the following estimate,
\begin{eqnarray}
&&-\dot{\xi}_{\ve}\int\frac{W'''(\theta_0)}{\varepsilon}(v_{\varepsilon}+\xi_{\ve})(\theta_0')^2\theta_1dy+
\dot{\xi}_{\ve}\int \frac{\Psi_0 }{\varepsilon}(\theta_0'(v_{\varepsilon}+\xi_{\ve}))'\theta_0'dy+\dot{\xi}_{\ve} \int \frac{Q_{\varepsilon}}{\varepsilon}(\theta_0')^2dy \nonumber\\
&&\hspace{40 pt}\leq \frac{C}{\delta\varepsilon^2}\int (\theta_0')^2(v'_{\varepsilon})^2dy+\frac{C}{\delta\varepsilon^2}\int \left\{\Qs^2+\Qc^2\right\}dy +\frac{C}{\delta}\xi^2_{\ve}+\delta\dot{\xi}_{\ve}^2.\label{simplf_beta_3}
\end{eqnarray}

\noindent
STEP 4. We use equalities \eqref{simplf_beta_0},\eqref{simplf_beta_1} and estimates \eqref{simplf_beta_2},\eqref{simplf_beta_3} in \eqref{beta1}, and take $\delta>0$ such that $\int (\theta_0')^2dy>4 \delta$ to derive \eqref{eq_for_dot_xi}.

\rightline{$\square$}

Now, multiplying \eqref{eq_for_dot_xi} by $|\beta|$ and adding to \eqref{cor_1} we obtain the following corrolary.
\begin{corrolary}
The following inequality holds for sufficiently small $|\beta|$
\begin{eqnarray}
\nonumber&&\frac{d}{2dt}\left[\int (\theta_0')^2 v_{\ve}^2 dy +c_0\xi^2_{\ve}+\frac{1}{\ve}\int \Qc^2 dy +\frac{1}{\ve}\int \Qs^2 dy \right]\\
\nonumber&&\hspace{40 pt}+\frac{1}{4\ve^2}\left[\int (\theta_0')^2(v'_{\ve})^2 dy +\int \Qc^2 dy +\int (\Qc')^2 dy +\int \Qs^2 dy +\int (\Qs')^2 dy\right]\\
\nonumber&&\hspace{60 pt}\leq c\xi^2_{\ve} +\left(c\beta^2-|{\beta}|\right)\dot{\xi}^2+\varepsilon^{2(N-1-\alpha)}\\
&&\hspace{80 pt}+ \int R_{\ve}\theta_0'(v_{\ve}+\xi_{\ve})dy+ |\dot{\xi}|\beta\int R_{\ve}\theta_0'dy. \label{cor_2}
\end{eqnarray}
\end{corrolary}

In the following lemma we obtain appropriate bounds for the last two terms in \eqref{cor_2}, i.e. terms 
containing $R_{\varepsilon}$.\\
\noindent{\bf Lemma 4.} {\it The following inequalities hold true% for $\alpha=3$ and $N\geq 4$:
\begin{list}{}{}
\item{(i)} for all $\alpha >2$ and $N\geq \alpha+1$ 
\begin{eqnarray}\nonumber
\int R_{\varepsilon}\theta_0'(v_{\varepsilon}+\xi_{\ve})dy&\leq & \frac{c}{\varepsilon} \int(\theta_0')^2 (v_{\varepsilon}')^2dy+c \xi_{\ve}^2+c\nonumber\\
&& +c\ve^{\alpha -2} \left(\int (\theta_0')^2 v_{\varepsilon}^2dy\right)^2+c\ve^{2\alpha -2}\left(\int (\theta_0')^2 v_{\varepsilon}^2dy\right)^3\nonumber\\
&& +c \ve^{\alpha -2} |\xi_{\ve}|^3+c \ve^{2\alpha -2}\xi^{4}_{\ve}\nonumber\\
&&+c\int \Qs^2 dy+c\int (\Qs')^2dy\nonumber\\
&&+c\int \Qc^2 dy+c\int (\Qc')^2dy. \label{r_eps_1}
\end{eqnarray}
\item{(ii)} for all $\alpha >2$ and $N\geq \alpha+1$
\begin{eqnarray}
\nonumber |\dot\xi_{\ve}| \int R_{\varepsilon}\theta_0'dy&\leq&   \frac{1}{8}\dot{\xi}_{\ve}^2+c\varepsilon^{4(\alpha-1)}\xi_{\ve}^{6}+c\left[\int (\theta_0')^2v_{\varepsilon}^2dy
+\xi_{\ve}^{2}\right]\\\nonumber&&+c\varepsilon^{2(\alpha-2)}\left[\int (\theta_0')^2v_{\varepsilon}^2dy+\xi_{\ve}^2\right]^2+c\int \Qs^2 dy+c\int \Qc^2 dy\\&&+
c\varepsilon^{2\alpha-2} \left(\int (\theta_0')^2v_{\varepsilon}^2dy\right)^{3/2} \left(\int (\theta_0')^2(v_{\varepsilon}')^2dy\right)^{1/2}|\dot{\xi}_{\ve}|\nonumber\\&&+\varepsilon^{\alpha-1}\left(\int u_{\varepsilon}^2dy \right)^{1/2}\left(\int Q_{\ve}^2+(Q'_{\ve})^2dy\right)^{1/2}|\dot{\xi}_{\ve}|. 
\label{r_eps_2}
\end{eqnarray}
%\item{(iii)} 
%\begin{equation}
%\int m_{\varepsilon}\Qs dy\leq c \int \Qs^2 dy +1  \label{q_eps}
%\end{equation}
\end{list} 
}
\noindent{P r o o f.}\\
To show (i) we estimate $\int R_{\ve}\theta_0'(v_{\ve}+\xi_{\ve})dy$ using the definition of $R_{\ve}$ \eqref{repr_for_R_eps}, the Cauchy-Schwarz inequality, Poincar\'e inequalities \eqref{poincare},\eqref{interp_3},\eqref{ii_4}, and estimates \eqref{general_Q} and \eqref{general_Q_2} from Lemma \ref{Rem_in_Q}, 
\begin{eqnarray*}
\int R_{\ve}\theta_0'(v_{\ve}+\xi_{\ve})dy &\leq& c \left[1+\ve^{2(N-1-\alpha)}+\ve^{\alpha}+\ve^{2\alpha-2}\right]\int (\theta_0')^2 (v'_{\ve})^2dy\\&&+c\left[1+\ve^{\alpha-1}+\ve^{2(N-\alpha)}\right]\xi_{\ve}^2+c\ve^{\alpha-2}|\xi_{\ve}|^3+c\ve^{2\alpha-2}\xi_{\ve}^4\\\
&&+c\ve^{\alpha -2}\left\{ \int (\theta_0')^2 v_{\ve}^2dy\right\}^2+c\ve^{2\alpha -2}\left\{\int (\theta_0')^2 v_{\ve}^2dy\right\}^3\\
&&+c\left[\int \Qs^2 dy +\int \Qc^2 dy \right]+c\ve^{\alpha-1}\left[\int (\Qs')^2 dy +\int (\Qc')^2 dy \right].
\end{eqnarray*}
Taking $\alpha >2$ and $N\geq \alpha+1$ we derive \eqref{r_eps_1}.

To prove  (ii) we note that definition  \eqref{repr_for_R_eps} of $R_{\ve}$, the Cauchy-Schwarz inequality, Poincar\'{e} inequalities \eqref{poincare} and \eqref{interp_3} yield
\begin{eqnarray}
 |\dot{\xi}_{\ve}|\int R_{\varepsilon}\theta_0' dy&\leq & \frac{1}{8}\dot \xi_{\ve}^2+ c (1+\varepsilon^{N-\alpha}){\xi}_{\ve}^2+c\ve^{4(\alpha-1)}\xi_{\ve}^6+c\int (\theta_0')^2 (v_{\varepsilon}')^2dy\nonumber\\
&&+c\varepsilon^{2(\alpha -2)}\left[\int (\theta_0')^2 v_{\ve}^2dy+\xi_{\ve}^2\right]^2\nonumber\\
&& +c\varepsilon^{2\alpha-2}\left(\int (\theta_0')^2v_{\ve}^2dy\right)^{3/2}\left(\int (\theta_0')^2(v_{\ve}')^2 dy \right)^{1/2}|\dot{\xi}_{\ve}|\nonumber\\
&&+c\int \Qs^2 dy+ c\int\Qc^2 dy\nonumber\\
&&+c\ve^{\alpha-1}\left(\int (\theta_0')^2 v_{\ve}^2dy +\xi_{\ve}^2\right)^{1/2}\left(\int Q_{\ve}^2dy +\int (Q_{\ve}')^2 dy\right)^{1/2}|\dot{\xi}|.  \nonumber
\end{eqnarray}

%\noindent The proof of (iii) is straightforward. \\
\rightline{$\square$}
Now it is convinient to introduce the folloowing notation, 
\begin{eqnarray*}
\mathcal{E}_{\varepsilon}(t)&=&\int (\theta_0')^2 v_{\varepsilon}^2 dy+c_0\xi_{\ve}^2+\frac{1}{\varepsilon}\int \Qs^2 dy+\frac{1}{\varepsilon}\int \Qc^2dy
\\
\mathcal{D}_{\varepsilon}(t)&=&\frac{1}{8\varepsilon^2}\left\{\int (\theta_0')^2(v_{\varepsilon}')^2 dy +\int\Qs^2dy +\int (\Qs')^2 dy+\int\Qc^2dy +\int (\Qc')^2 dy \right\}\\&&+\left(\frac{|\beta|}{8}-c\beta^2\right)\dot{\xi}_{\ve}^2.
\end{eqnarray*}
In terms of $\mathcal{E}_{\varepsilon}$ and $\mathcal{D}_{\varepsilon}$ we can rewrite \eqref{r_eps_2} for $\alpha =3$ and $N=4$ 
 in the following form, 
\begin{equation*}
\int R_{\varepsilon}\theta_0'dy\dot{\xi}_{\ve}\leq c\mathcal{E}_{\varepsilon}+c\varepsilon\mathcal{E}^2_{\varepsilon}+c\varepsilon^6\mathcal{E}^3_{\varepsilon}+c(\varepsilon+\varepsilon^2\mathcal{E}^{1/2}_{\varepsilon}+\varepsilon^{4}\mathcal{E}^{3/2})
\mathcal{D}_{\varepsilon}
\end{equation*}
Substituting the above inequality and  \eqref{r_eps_2} into \eqref{cor_2} we obtain 
\begin{equation}\nonumber
\dot{\mathcal{E}}_{\varepsilon}+\frac{1}{2}\mathcal{D}_{\varepsilon}\leq c\mathcal{E}_{\varepsilon}+c\varepsilon \mathcal{E}^{3/2}_{\varepsilon}+c\varepsilon\mathcal{E}^2_{\varepsilon}+c\varepsilon^6\mathcal{E}^3_{\varepsilon}
+c(\varepsilon+\varepsilon^2\mathcal{E}^{1/2}_{\varepsilon}+\varepsilon^{4}\mathcal{E}^{3/2}_{\varepsilon})
\mathcal{D}_{\varepsilon}.
\end{equation}
Assume that  
\begin{equation}\nonumber
\mathcal{E}_{\varepsilon}(t)\leq c\text{ for all }t\in [0,t^{\star}).
\end{equation}
Then for $t\in[0,t^{\star})$ we have
\begin{equation}
\dot{\mathcal{E}}_{\varepsilon}\leq c \mathcal{E}_{\varepsilon}
\end{equation}
Thus,  $[0,T]\subset [0,t_{\star})$ sufficiently small $\varepsilon$.
 {\bf This concludes the proof of Theorem \ref{theorem_sharp}.}

%Assume that 
%\begin{equation}\nonumber
%\mathcal{E}_{\varepsilon}(t)\leq c_1\text{ for all }t\in [0,t^{\star}).
%\end{equation}
%Then for $t\in[0,t^{\star})$ we have
%\begin{equation}
%\dot{\mathcal{E}}_{\varepsilon}\leq 3C^2 \mathcal{E}_{\varepsilon}
%\end{equation}
%Thus, for $t\in[0,t^{\star})$
%\begin{eqnarray}
%\mathcal{E}_{\varepsilon}(t)&\leq &\mathcal{E}_{\varepsilon}(0)e^{3C^2t}
%\end{eqnarray}
%Thus, we need to consider $\varepsilon$ so small such that $t^{\star}>T$. {\bf This concludes the proof of a priori bounds.}
%\newpage

%%%%%%%%%%%%%%%%%%%%%%%%%%%%%%%%%%%%%%%%%%%%%%%%%%%%
%%%%%%%%%%%% Section: Formal derivation of curvatur motion %%%%%%%%%%%%%%%%%%%%
%%%%%%%%%%%%%%%%%%%%%%%%%%%%%%%%%%%%%%%%%%%%%%%%%%%%
\section{Formal derivation of \eqref{motion1}}\label{formalderivation}

In this section we formally derive equation \eqref{motion1} for 2D   system  (\ref{eq1}-\ref{eq2}) with gradient coupling.  Analogous  derivation for the single  Allen-Cahn equation has been done in  \cite{RubSterKel89}, %\textcolor{red} {\cite{MatSha}} 
and  \cite{CheHilLog10}. 
% assuming that its rigorous justification is similar to 1D done in the Section \ref{limit}. Following \cite{RubSterKel89} we introduce local coordinates. These coordinates are convenient in this problem because main variations of the functions $\rho_{\ve}$ and $P_{\ve}$ are in the direction perpenicular to the interface of the scale $\ve$. 

Consider a subdomain $\omega_t \subset \Omega \subset \mathbb{R}^2$ ($\omega_t$ is occupied by the cell) so that $\Gamma(0)=\partial \omega_0$ (boundary of the cell at $t=0$). For all $t\in[0,T]$ consider a closed curve $\Gamma(t)$ s.t. $\partial \omega _t=\Gamma(t)$ and $\omega_t\subset \Omega$. Let $X_0(s,t)$ be a parametrization of $\Gamma(t)$. In a vicinity of $\Gamma(t)$  the parameters $s$ and the  signed distance   $r$ to $\Gamma(t)$  will be used as local coordinates, so that
\begin{equation}\nonumber
x=X_0(s,t)+r \n(s,t)=X(r,s,t), \qquad \text{where  $\n$ is an inward normal.}
\end{equation}
The inverse mapping to $x=X(r,s,t)$ is given by 
\begin{equation}\nonumber
r=\pm\text{dist}(x,\Gamma(t)),\;\;s=S(x,t),
\end{equation}
where in the formula for $r$ we choose $-$ if $x\in \omega_t$ and $+$, if $x\notin \omega_t$. 
Recall that $\Gamma(t)$ is the limiting location of interface as $\ve \to 0$. For fixed $\ve$ we are looking for interface in the form of $\ve$-perturbation of $\Gamma(t)$: 
\begin{equation}\nonumber
\tilde{\Gamma}_{\ve}(t)=\Gamma(t)+\varepsilon h_{\varepsilon}(s,t).
\end{equation}
Introduce the limiting velocity $V_0:=-\partial_tr$ and the distance to $\Gamma_{\ve}(t)$ rescaled by $\ve$:
 %Introduce rescaled  variable: 
\begin{equation}\label{def_of_z}
z=z^{\varepsilon}(x,t)=\frac{r-\varepsilon h_{\varepsilon}(S(x,t),t)}{\varepsilon}.
\end{equation}
Next we define a rule that  for all positive $t$ transforms any given function $w$ of original variable $x$ into the corresponding function $\tilde{w}$ in local coordinates $(z,s)$:
\begin{equation}\nonumber
w(x,t)=\tilde{w}\left(\frac{r(x,t)-\varepsilon h_{\varepsilon}(S(x,t),t)}{\varepsilon},S(x,t),t\right).
\end{equation} 
By applying this rule for the functions $\rho_{\ve}$ and $P_{\ve}$ we define $\tilde{\rho}_{\ve}$ and $\tilde{P}_{\ve}$:
\begin{equation}\nonumber
\tilde{\rho}_{\varepsilon}(z,s,t)=\rho_{\varepsilon}(x,t)\text{ and }\tilde{P}_{\varepsilon}(z,s,t)=P_{\varepsilon}(x,t).
\end{equation}
We now introduce asymptotic expansions in local coordinates:
\begin{eqnarray} \label{expansion_5}
\tilde\rho_{\varepsilon}(z,s,t)&=&\theta_0(z,s,t)+\varepsilon \theta_{1}(z,s,t)+...\\
\tilde P_{\varepsilon}(z,s,t)&=&\Psi_0(z,s,t)+\varepsilon \Psi_{1}(z,s,t)+...\\
h_{\varepsilon}(s,t)&=& h_1(s,t)+\varepsilon h_2(s,t)+...\\
\lambda_{\varepsilon}(t)&=&\frac{\lambda_0(t)}{\ve}+\lambda_1(t)+\ve\lambda_2(t)+...\label{expansion9}
\end{eqnarray}
Now, substitute (\ref{expansion_5}-\ref{expansion9}) into \eqref{eq1} and \eqref{eq2}. Equating coefficients at $\ve^{-2},\ve^{-1}$ and ${\ve}^0$, we get:
\begin{equation}\label{theta_0z}
\frac{\partial^2 \theta_0}{\partial z^2}=W'(\theta_0),
\end{equation}
and 
\begin{eqnarray}
-V_0\frac{\partial\theta_0}{\partial z}&=&\frac{\partial^2 \theta_1}{\partial z^2}-W''(\theta_0)\theta_1+\frac{\partial \theta_0}{\partial z}\kappa(s,t)+(\Psi_0\cdot \n)\frac{\partial \theta_{0}}{\partial z}+\lambda_{0}(t),\label{solv_for_theta_1}\\
-V_0\frac{\partial \Psi_0}{\partial z}&=& \frac{\partial^2 \Psi_0}{\partial z^2}-\Psi_0+\beta \frac{\partial \theta_0}{\partial z}\n.\label{def_of_psi}
\end{eqnarray}
where $\kappa(s,t)$ is the curvature of $\Gamma_0(t)$. The curvature $\kappa$ appears in the equation when one rewrites the Laplace operator in \eqref{eq1} in local coordinates $(z,s)$. The solvability condition for the equation for $\theta_1$ \eqref{solv_for_theta_1} %(orthogonality to $\theta_0^\prime$) turns into
yields
\begin{equation}\label{formalmotion1}
c_0 V_0=c_0 \kappa(s,t)+\int (\Psi_0\cdot \n)\left(\frac{\partial \theta_0}{\partial z}\right)^2 dz
+ \lambda_0(t).
\end{equation} 
From the definition \eqref{lagrange} it follows that $\int_{\Omega}\partial_t \rho_{\ve}=0$. Substitute expansions for $\rho_{\ve}$ \eqref{expansion_5} into $\int_{\Omega}\partial_t \rho_{\ve}=0$, take into account that $\theta_0$ does not explicitly depend on $s$ and $t$ (which follows from the equation \eqref{theta_0z}; note that $\theta_0$ still depends on $t$ implicitly, through variable $z$ which by \eqref{def_of_z} is a function of $t$). By integrating RHS of \eqref{formalmotion1} we get 
 
\begin{equation}\label{def_lambda_0}
\lambda_0(t)=-\int \left\{c_0 \kappa(s,t)+\int (\Psi_0\cdot \n)\left(\frac{\partial \theta_0}{\partial z}\right)^2 dz\right\} ds. 
\end{equation}  

Introduce $\tilde{\Psi}_0:=\Psi_0/\beta$. Since equation \eqref{def_of_psi} is linear with respect to $\Psi_0$ and the inhomogenuity is linearly proportional to $\beta$,  function $\tilde{\Psi}_0$ does not depend on $\beta$.  Finally, define 
\begin{equation}\label{def_for_motion1}
\Phi (V) : = \int (\tilde{\Psi}_0\cdot \n)\left(\frac{\partial \theta_0}{\partial z}\right)^2dz.
\end{equation}

By substituting \eqref{def_for_motion1} and \eqref{def_lambda_0} into equation \eqref{formalmotion1} we derive \eqref{motion1}.

%%%%%%%%%%%%%%%%%%%%%%%%%%%%%%%%%%%%%%%%%%%%%%%%%%%%
%%%%%%%%%%%%   Appendix   %%%%%%%%%%%%%%%%%%%%%%%%%%%%%%%%%%
%%%%%%%%%%%%%%%%%%%%%%%%%%%%%%%%%%%%%%%%%%%%%%%%%%%%
\appendix
\section {Appendix}
\subsection {Derivation of equations for $u_{\varepsilon}$ and $Q_{\varepsilon}.$} %{\eqref{eq_for_u}, \eqref{eq_for_Q}}.}
% and definition of functions  $a_\varepsilon(t,y),$ $b_{k,\varepsilon}(t,y),k=1,2,3$, $e_{\varepsilon}(t,y),g_{\varepsilon}(t,y)$, $m_{\varepsilon}(t,y)$.}
First, rewrite the linear part of the equation: 
\begin{eqnarray}
&&\frac{\partial \rho_{\varepsilon}}{\partial t}-\partial^2_{x}\rho_{\varepsilon}-\frac{F(t)}{\varepsilon}=\varepsilon^{\alpha}\left\{\frac{\partial u_{\varepsilon}}{\partial t}+\frac{V_0u_{\varepsilon}'}{\varepsilon}-\frac{u_{\varepsilon}''}{\varepsilon^2}\right\}+\varepsilon^{\alpha}\frac{V-V_0}{\varepsilon}u_{\varepsilon}\nonumber\\
&&\hspace{60 pt}-\frac{\theta_0''}{\varepsilon^2}+\frac{-\theta_1''+V_0\theta_0'-F_0}{\varepsilon}+...+\varepsilon^{N-2}\left(\dot\theta_{N-2}-\theta_N''+\sum\limits_{j=0}^{N-1}V_j\theta_{N-1-j}'\right)\nonumber\\
&&\hspace{60 pt}+\varepsilon^{N-1}\left(-\dot{\theta}_{N-1}+\sum\limits_{j=0}^{N}V_j\theta_{N-j}'\right)+\varepsilon^N r_1(\varepsilon,y,t)
%\\
%&&\hspace{60 pt}+\sum\limits_{i=0}^{N}\left(V-\sum_{j=0}^{N-i}\varepsilon^iV_i\right)\varepsilon^{i-1}\theta_i'-\sum\limits_{i=N}^{\infty}\varepsilon^{i-1}F_i.
\label{for_linear}
\end{eqnarray} 
Here 
\begin{eqnarray}\nonumber
\varepsilon^{N} r_1(\varepsilon,y,t)= \varepsilon^{N}\left[\dot{\theta}_N+\sum\limits_{k=0}^{N}\varepsilon^{k}\left\{\sum\limits_{j=1+k}^{N+1+k}V_j\theta_{N+1+k-j}'\right\}\right]
\end{eqnarray}
To analyze nonlinear part denote 
\begin{equation}
h_{\varepsilon}(t,y)=\sum\limits_{i=1}^{N}\varepsilon^{i-1}\theta_i(t,y)\text{ and }r_{\varepsilon}(t,y)=\sum_{i=1}^{N}\varepsilon^{i-1}\Psi_{i}(t,y).
\end{equation}
Thus, $\rho_{\varepsilon}=\theta_0+\varepsilon h_{\varepsilon}+\varepsilon^{\alpha} u_{\varepsilon}$ and $P_{\varepsilon}=\Psi_0+\varepsilon r_{\varepsilon}+\varepsilon^{\alpha} Q_{\varepsilon}$. Rewrite nonlinear terms:
\begin{eqnarray}
P_{\varepsilon}\partial_{x}\rho_{\varepsilon}&=&\varepsilon^{-1}(\Psi_0+\varepsilon r_{\varepsilon}+\varepsilon^{\alpha}Q_{\varepsilon})(\theta_0'+\varepsilon h_{\varepsilon}'+\varepsilon^{\alpha}u_{\varepsilon}')\nonumber\\
&=&\varepsilon^{2\alpha-1}Q_{\varepsilon}u'_{\varepsilon}+\varepsilon^{\alpha-1}\Psi_0u'_{\varepsilon}+
\varepsilon^{\alpha-1}\theta_0'Q_{\varepsilon}+\varepsilon^{\alpha}r_{\varepsilon}u'_{\varepsilon}
+\varepsilon^{\alpha}h_{\varepsilon} Q_{\varepsilon}\nonumber\\
&&+\sum\limits_{i=1}^{N}\varepsilon^{i-2}\sum\limits_{j=0}^{i-1}\Psi_{j}\theta_{i-1-j}'+\varepsilon^{N-1}\sum\limits_{j=0}^{N}\Psi_j\theta_{N-j}'+\varepsilon^{N}r_2({\varepsilon},y,t).\label{for_P_rho}
%\sum\limits_{i=N+1}^{2N}\varepsilon^{i-N-1}\sum\limits_{j=i-N}^{N}\Psi_j\theta_{i-j}'
\end{eqnarray}
Here  
\begin{equation}\nonumber
\varepsilon^{N}r_2(\varepsilon,y,t)=\varepsilon^{N}\sum\limits_{i=N+1}^{2N}\varepsilon^{i-1-N}\sum\limits_{j=i-N}^{N}\Psi_j\theta_{i-j}'.
\end{equation}
\begin{eqnarray}
\nonumber\frac{W'(\rho_{\varepsilon})}{\varepsilon^2}&=&
\left\{W''(\theta_0)+\varepsilon W'''(\theta_0)h_{\varepsilon}+\varepsilon^2\frac{W^{\text{(iv)}}(\theta_0)}{2}h_{\varepsilon}^2\right\}\varepsilon^{\alpha-2}u_{\varepsilon}\\\nonumber
&&+\left\{\frac{W'''(\theta_0)}{2}+\varepsilon\frac{W^{\text{(iv)}}(\theta_0)}{2}h_{\varepsilon}\right\}\varepsilon^{2\alpha-2}u_{\varepsilon}^2+\frac{W^{\text{(iv)}}(\theta_0)}{6}\varepsilon^{3\alpha-2}u_{\varepsilon}^3.\\\nonumber
&&+\frac{W'(\theta_0)}{\varepsilon^2}+\frac{W''(\theta_0)\theta_1}{\varepsilon}+...+\varepsilon^{N-2}\left[W''(\theta_0)\theta_N+(dW)^{(N)}\right]\\
&&+\varepsilon^{N-1}(dW)^{(N+1)}\label{for_W}+\varepsilon^{N} r_3(\varepsilon,y,t),
\end{eqnarray}
Here 
\begin{equation}\nonumber
\varepsilon^{N}r_3(\varepsilon,y,t)=\varepsilon^{N}\sum\limits_{k=2}^{2N}\varepsilon^{k-2}
\left\{\sum\limits_{\begin{array}{c}i_1+i_2=N+k\\1\leq i_1,i_2\leq N\end{array}}\frac{W'''(\theta_0)}{2}\theta_{i_1}\theta_{i_2}+
\sum\limits_{\begin{array}{c}i_1+i_2+i_3=N+k\\1\leq i_1,i_2,i_3\leq N\end{array}}\frac{W^{(iv)}(\theta_0)}{6}\theta_{i_1}\theta_{i_2}\theta_{i_3}
\right\}
\end{equation}
Summing together \eqref{for_linear},\eqref{for_P_rho}, \eqref{for_W} and dividing by $\varepsilon^{\alpha}$ we get \eqref{eq_for_u}. 
Derivation \eqref{eq_for_Q} is simple since \eqref{orig_2} is linear.

In the end of this subsection we write expression for functions in \eqref{repr_for_R_eps}. 
\begin{equation}\nonumber
a_{\varepsilon}(t,y)=-\dot{\theta}_{N-1}+\sum\limits_{j=0}^{N}V_j\theta_{N-j}'-F_N+\sum\limits_{j=0}^{N}\Psi_{j}\theta_{N-j}'-(dW)^{(N+1)},
\end{equation}
\begin{equation}\nonumber
b_{0,\varepsilon}(t,y)=-r_1(\ve,y,t)+r_2(\ve,y,t)-r_3(\ve,y,t),
\end{equation}
\begin{equation}\nonumber
b_{1,\ve}(t,y)=-W'''(\theta_0)(\theta_2+\varepsilon\theta_3+...+\varepsilon^{N-2}\theta_N)+\frac{W^{\text{(iv)}}(\theta_0)}{2}h_{\ve}^2,
\end{equation}
\begin{equation}\nonumber
b_{2,\ve}(t,y)=\frac{W'''(\theta_0)}{2}+\varepsilon\frac{W^{\text{(iv)}}(\theta_0)}{2}h_{\ve},
\end{equation}
\begin{equation}\nonumber
b_{3,\ve}(t,y)=\frac{W^{\text{(iv)}}(\theta_0)}{6},
\end{equation}
\begin{equation}\nonumber
e_{\ve}(t,y)=\frac{V-V_0}{\varepsilon}+r_{\ve},
\end{equation}
\begin{equation}\nonumber
g_{\varepsilon}(t,y)=h'_{\ve}.
\end{equation}
and the function $m_{\ve}$ from \eqref{eq_for_Q}:
\begin{equation}\nonumber
m_{\varepsilon}(t,y)=\sum\limits_{k=0}^{N}\varepsilon^{k}\left[\sum\limits_{j=k+1}^{N+k+1}V_j\Psi_{N+k+1-j}\right]
\end{equation}

\subsection{Auxiliary Inequalities.}
\noindent{\bf Assumption on $(\theta_0')^2$.} There exist $\kappa>0$ and $c_0>1$ such that 
\begin{equation}\label{exp}
c_0^{-1}e^{-\kappa |y|}<(\theta_0'(y))^2\leq c_0 e^{-\kappa |y|},\;\;\;y\in \mathbb R
\end{equation}

\medskip

\noindent{\bf Remarks.} 1. All inequalities below will be proven for particular case $(\theta_0'(y))^2=e^{-\kappa y}$. The result is obviously extended for all $(\theta'_0)^2$ satisfying \eqref{exp}.\\
% We will use below the following inequality for $x>0$:
%\begin{equation}\label{aux_theta}
%\frac{c_0^{-2}}{\kappa}(\theta_0'(x))^2\leq\int_{x}^{\infty}(\theta_0'(t))^2dt\leq \frac{c_0^2}{\kappa} (\theta_0'(x))^2.
%\end{equation}
\noindent 2. It is easy to make sure that $\theta_0'$ from the cell movement problem satisfies \eqref{exp}. Indeed, 
\begin{equation}
\theta_0(y)=\frac{1}{2}(1-\tanh \frac{y}{2\sqrt{2}})\text{ and }(\theta'_0(y))^2=\frac{1}{32}\frac{1}{\cosh^4\frac{y}{2\sqrt{2}}}
\end{equation}
and $2e^{-|x|}\leq \cosh x\leq 4 e^{-|x|}$ for $x\in \mathbb R$.

\medskip 

\noindent{\bf Theorem.} ({\it Poincar\'e inequality})
\begin{equation}\label{poincare}
\int (\theta_0')^2 (v-<v>)^2dy\leq c_P\int (\theta_0')^2(v')^2dy, \;\;v\in C^1(\mathbb R)\cap L^{\infty}(\mathbb R),
\end{equation}
where 
\begin{equation}
<v>=\frac{\int(\theta_0')^2vdy}{\int(\theta_0')^2dy}.
\end{equation}
{P r o o f.}\\
\noindent{STEP 1.} ({\it Friedrich's inequality.}) Take any $u\in C^1(\mathbb R)\cap L^{\infty}(\mathbb R)$ s.t. $u(0)=0$ and let us prove the Friedrich's inequality:
\begin{equation}\label{friedrich}
\int (\theta_0')^2u^2dy\leq c_F\int(\theta_0')^2(u')^2dy,
\end{equation} 
where $c_F$ does not depend on $u$.
Indeed,
\begin{eqnarray}
\int_0^{\infty} e^{-\kappa y} u^2 dy &=& \int_0^{\infty} \left(\int_{y}^{\infty}e^{-\kappa t}dt\right) u'u dy\nonumber\\ &\leq& \int_{0}^{\infty}\left(\int_{y}^{\infty}e^{-\kappa t}dt\right) |u'||u|dy\nonumber\\
&=&\frac{1}{\kappa}\int_0^{\infty}e^{-\kappa y}|u'||u|dy\nonumber\\
&\leq & \frac{1}{\kappa}\left(\int_{0}^{\infty} e^{-\kappa y}(u')^2 dy\right)^{1/2}\left(\int_{0}^{\infty} e^{-\kappa y}u^2 dy\right)^{1/2}.\nonumber
\end{eqnarray}
Thus, 
\begin{equation}
\int_0^{\infty} (\theta_0')^2 u^2 dy\leq \frac{c_0^2}{\kappa}\left(\int_{0}^{\infty} (\theta_0')^2(u')^2 dy\right)^{1/2}\left(\int_{0}^{\infty} (\theta_0')^2u^2 dy\right)^{1/2}.\nonumber
\end{equation}
Dividing the latter inequality by $\left(\int_{0}^{\infty} (\theta_0')^2u^2 dy\right)^{1/2}$, and than taking square of both sides we have  
\begin{equation}\nonumber
\int_{0}^{\infty}(\theta_0')^2u^2dy \leq \frac{c_0^4}{\kappa^2}\int_0^{\infty}(\theta_0')^2(u')^2 dy.
\end{equation}
Of course, similar inequality is valid on $(-\infty,0)$: 
\begin{equation}\nonumber
\int_{-\infty}^{0}(\theta_0')^2u^2dy \leq \frac{c_0^4}{\kappa^2}\int_{-\infty}^{0}(\theta_0')^2(u')^2 dy.
\end{equation}
Hence we proven the Friedrich's inequality \eqref{friedrich} with the constant $c_F=c_0^4/\kappa^2$. \\
\noindent{STEP 2.} We prove the Poincar\'e inequality \eqref{poincare} by contradiction: assume that there exists a sequence $v_n\in C^1(\mathbb R)\cap L^{\infty}(\mathbb R)$ such that 
\begin{equation}\nonumber
\int (\theta_0')^2v^2_ndy =1, \;\;\int(\theta_0')^2 v_n dy=0\text{ and }\int (\theta_0')^2 (v_n')^2dy\rightarrow 0.  
\end{equation}
Apply Friedrich's inequality \eqref{friedrich} for functions $v_n(y)-v_n(0)$: 
\begin{equation}\nonumber
\int (\theta_0')^2(v_n(y)-v_n(0))^2 dy \leq c_F\int (\theta_0')^2 (v_n')^2 dy\rightarrow 0.
\end{equation} 
On the other hand, 
\begin{equation}\nonumber
\int (\theta_0')^2(v_n(y)-v_n(0))^2 dy=\int (\theta_0')^2 v_n^2dy +v_n^2(0)\int(\theta_0')^2 dy\geq \int (\theta_0')^2 v_n^2dy.  
\end{equation}
Hence, 
\begin{equation}\nonumber
\int (\theta_0')^2 v_n^2 dy \rightarrow 0
\end{equation}
which contradicts to $\int (\theta_0')^2 v_n^2 dy=1$. The contradiction proves the theorem.\\
\rightline {$\square$}

\medskip

\noindent{\bf Theorem.} ({\it Interpolation inequality})
\begin{equation}\label{interp_3}
\int (\theta_0')^3 v^3dy\leq c_{I} \left(\int (\theta_0')^2\left\{v^2+(v')^2\right\} dy\right)^{1/2}\left(\int (\theta_0')^2 v^2 dy\right)
\end{equation}
{P r o o f.}\\
Consider first $y>0$. 
\begin{eqnarray}
e^{-\kappa y} |v(y)|&=&|\kappa\int_y^{\infty}e^{-\kappa x} v(x)dx-\int_y^{\infty}e^{-\kappa x} v'(x)dx|\nonumber\\&\leq &
\left(\int_0^{\infty}e^{-\kappa y}v^2dy\right)^{1/2}\left(\int_{y}^{\infty}e^{-\kappa y} dy\right)^{1/2}\nonumber\\&&+ 
\frac{c_0}{\kappa^{1/2}}\left(\int_0^{\infty}e^{-\kappa y}(v')^2 dy\right)^{1/2}\left(\int_{y}^{\infty}e^{-\kappa y}dy\right)^{1/2}\nonumber\\
&\leq &c_1\left(\int_{0}^{\infty} (\theta_0')^2\left\{v^2+(v')^2 \right\}dy \right)^{1/2}|e^{-\kappa y/2}| \nonumber
\end{eqnarray}
Thus, we have 
\begin{equation}\label{last_step}
e^{-\kappa y}|v|\leq c_2\left(\int_{0}^{\infty} e^{-\kappa t}\left\{v^2+(v')^2 \right\}dt \right)^{1/2}e^{-\kappa y/2}
\end{equation}
Multiply by $e^{-\kappa y/2}|v|^2$ and intergrate over $(0,\infty)$: 
\begin{equation}\nonumber
\int_{0}^{\infty}e^{-3\kappa y/2}|v|^3 dy \leq c_2 \left(\int e^{-\kappa t}\left\{v^2+(v')^2 \right\}dy\right)^{1/2}\int e^{-\kappa y}v^{2}dy.
\end{equation}
Rederiving the same estimate for $(-\infty,0)$ we prove the theorem.\\\rightline {$\square$} 
%\newpage
%\medskip 

%\noindent{\bf Remark.} Take $n\geq 3$ and in the last step of the proof of the theorem let us multiply \eqref{last_step} by $|\theta_0'|^{n-2}|v|^{n-1}$:
%\begin{eqnarray}
%\int |\theta_0'|^n|v|^ndy&\leq& c_2\left(\int (\theta_0')^2\left\{v^2+(v')^2 \right\}dy\right)^{1/2}\int (\theta_0')^{n-1}v^{n-1}dy\\
%&&...\nonumber\\
%&\leq &c_2^{n-2}\left(\int (\theta_0')^2\left\{v^2+(v')^2 \right\}dy\right)^{(n-2)/2}\int (\theta_0')^2v^{2}dy
%\end{eqnarray}
%Moreover if $<v>=0$ we can use Poincare inequality and thus: 
%\begin{eqnarray}
%\int |\theta_0'|^n|v|^ndy \leq c_n\left(\int(\theta_0')^2 (v')^2 dy\right)^{(n-2)/2}\int (\theta_0')^2 v^2 dy 
%\end{eqnarray}
%However, we need another interpolation inequality with power 1/2 of $H^1$-norm. 
%
%\medskip 
%

\noindent{\bf Theorem.} ({\it Interpolation inequality for $n=4$}):
\begin{equation}\label{ii_4}
\int |\theta_0'|^4|v|^4dy \leq c_4\left(\int(\theta_0')^2 \left\{v^2+(v')^2\right\} dy\right)^{1/2}\left(\int (\theta_0')^2 v^2 dy\right)^{3/2} 
\end{equation}
\noindent{P r o o f.}\\
%\noindent STEP 1. ({\it the inequality for exponential weight})%Prove the inequality \eqref{ii_n} for exponential weight.
\noindent Take $y>0$
\begin{eqnarray*}
e^{-2\kappa y}v^{2}(y) &=& 2\kappa\int_{y}^{\infty}e^{-2\kappa t}v^{2}dt - 2 \int_{y}^{\infty}e^{-2\kappa t}vv'dt\\
&\leq &2\kappa \int _0^{\infty}e^{-2\kappa t}v^2 dt +2 \left(\int_0^{\infty} e^{-2\kappa t}(v')^2dt\right)^{1/2}\left(\int_{0}^{\infty}e^{-2\kappa t}v^2dt\right)^{1/2} 
\\ &\leq & c\left(\int_{0}^{\infty}e^{-2\kappa t} \left\{v^2+(v')^2\right\}dt\right)^{1/2}\left(\int_0^{\infty}e^{-2\kappa t}v^2dt\right)^{1/2}
\end{eqnarray*}
Multiply by $e^{-2\kappa y}v^2(y)$ and integrate over positive $y$:
\begin{equation*}
\int_{0}^{\infty} e^{-4\kappa y}v^4(y)dy\leq  c\left(\int_{0}^{\infty}e^{-2\kappa t} \left\{v^2+(v')^2\right\}dt\right)^{1/2}\left(\int_0^{\infty}e^{-2\kappa t}v^2dt\right)^{3/2}.
\end{equation*}
Use \eqref{exp} and the same inequality for $(-\infty,0)$. \\\rightline{$\square$}
%It is possible to prove that 
%\begin{equation}\label{l6}
%\int (\theta_0')^6v_{\varepsilon}^6dy \leq C \left(\int (\theta_0')^2(v_{\varepsilon}')^2dy\right)^2+\left(\int (\theta_0')^2v_{\varepsilon}^2dy\right)^4
%\end{equation}

\subsection{Uniqueness of original problem \eqref{eq1}-\eqref{eq2}}
In this appendix we prove uniqueness of the solution to the problem \eqref{eq1}-\eqref{eq2} in the following class: 
\begin{equation}\label{class_for_rho}
\rho_{\varepsilon}^{(i)}\in C([0,T]\times \Omega)\cap C([0,T];H^1(\Omega)\cap L^4(\Omega)),\;\partial_{t}\rho_{\varepsilon}^{(i)}\in L^2((0,T)\times \Omega).
\end{equation}
\begin{equation}\label{class_for_P}
P_{\varepsilon}^{(i)}\in C([0,T];  L^4(\Omega))\cap L^2(0,T;H^1(\Omega)).
\end{equation}

Assume that $\rho_{\varepsilon}^{(1)}, P_{\varepsilon}^{(1)}$ and $\rho_{\varepsilon}^{(2)}, P_{\varepsilon}^{(2)}$ are solutions satisfying \eqref{class_for_rho} and \eqref{class_for_P} for some positive $T>0$.

\noindent Take $\overline{\rho}_{\varepsilon}=\rho_{\varepsilon}^1-\rho_{\varepsilon}^2$ and ${{\overline{P}_{\varepsilon}}}=P_{\varepsilon}^1-P_{\varepsilon}^2$. Then equation for ${{\overline{\rho}_{\varepsilon}}}$:
\begin{equation}\nonumber
\frac{\partial {{\overline{\rho}_{\varepsilon}}}}{\partial t}=\Delta {{\overline{\rho}_{\varepsilon}}} +a(x,t){{\overline{\rho}_{\varepsilon}}}-b(x,t)\cdot \nabla {{\overline{\rho}_{\varepsilon}}}+c(x,t)\cdot {{\overline{P}_{\varepsilon}}}-\overline{\lambda}(t),
\end{equation} 
where 
\begin{equation}\nonumber
a(x,t)= \int\limits_0^1 W''(\rho_{\varepsilon}^{(1)}+s{\overline{\rho}_{\varepsilon}})ds,\;\;
b(x,t) = P_{\varepsilon}^{(1)},\;\;
c(x,t)=\nabla \rho_{\varepsilon}^{(2)},\;\;
\overline{\lambda}(t)= \lambda _1(t)-\lambda_2(t).
\end{equation}
We know that $|a(x,t)|<c$, $\|b(\cdot,t)\|_{H^1(\Omega)}+\|b(\cdot,t)\|_{L^4(\Omega)}<c$, $\|c(\cdot,t)\|_{L^2(\Omega)}\leq c$ and 
\begin{equation}\nonumber
\overline{\lambda}(t)=\frac{1}{|\Omega|}\int_{\Omega}\left\{a(x,t){{\overline{\rho}_{\varepsilon}}}-b(x,t)\cdot \nabla {{\overline{\rho}_{\varepsilon}}}+c(x,t)\cdot {{\overline{P}_{\varepsilon}}}\right\}dx
\end{equation}
and thus,
\begin{equation}\nonumber
|\overline{\lambda}(t)|\leq c\|{{\overline{\rho}_{\varepsilon}}}\|+c\|\nabla {{\overline{\rho}_{\varepsilon}}}\|+c\|{{\overline{P}_{\varepsilon}}}\|
\end{equation}
And energy estimate is (if multiply by ${{\overline{\rho}_{\varepsilon}}}$):
\begin{equation}\nonumber
\frac{\text{d}}{\text{d}t}\left[\int{{\overline{\rho}_{\varepsilon}}}^2\right]+2\int |\nabla {{\overline{\rho}_{\varepsilon}}}|^2 \leq c\| {{\overline{\rho}_{\varepsilon}}}\|
^2+c\|b\|_{L^4}\|\nabla {{\overline{\rho}_{\varepsilon}}}\|_{L^2}\|{{\overline{\rho}_{\varepsilon}}}\|_{L^4}+c\|c\|_{L^2}\|{{\overline{P}_{\varepsilon}}}\|_{L^4}\|{{\overline{\rho}_{\varepsilon}}}\|_{L^4}+c|\overline{\lambda}(t)|\|{{\overline{\rho}_{\varepsilon}}}\|
\end{equation}
Thus, using interpolation inequality \begin{equation}\nonumber\|\nabla {{\overline{\rho}_{\varepsilon}}}\|\|{{\overline{\rho}_{\varepsilon}}}\|_{L^4}\leq c\|\nabla {{\overline{\rho}_{\varepsilon}}}\|^{3/2}\|{{\overline{\rho}_{\varepsilon}}}\|^{1/2}\leq \nu \|\nabla {{\overline{\rho}_{\varepsilon}}}\|^2+c\|{{\overline{\rho}_{\varepsilon}}}\|^2.\end{equation}
\begin{equation}\nonumber
\|{{\overline{\rho}_{\varepsilon}}}\|_{L^4}\leq c\|{{\overline{\rho}_{\varepsilon}}}\|^{1/2}\|\nabla {{\overline{\rho}_{\varepsilon}}}\|^{1/2}<\nu \|\nabla {{\overline{\rho}_{\varepsilon}}}\|+c\|{{\overline{\rho}_{\varepsilon}}}\|
\end{equation}
and 
\begin{equation}\nonumber
|\overline{\lambda}(t)|\|{{\overline{\rho}_{\varepsilon}}}\|\leq \nu \|\nabla {{\overline{\rho}_{\varepsilon}}}\|^2+c\|{{\overline{\rho}_{\varepsilon}}}\|^2+c\|{{\overline{P}_{\varepsilon}}}\|^2
\end{equation}
we get
\begin{equation}\label{en_1}
\frac{\text{d}}{\text{d}t}\left[\int{{\overline{\rho}_{\varepsilon}}}^2\right]+\int |\nabla {{\overline{\rho}_{\varepsilon}}}|^2 \leq c\|{{\overline{\rho}_{\varepsilon}}}\|^2+c\|{{\overline{P}_{\varepsilon}}}\|_{L^4}^2.
\end{equation}
%\begin{equation}\label{en_1}
%\frac{\text{d}}{\text{d}t}\left[\int(\nabla {{\overline{\rho}_{\varepsilon}}})^2\right]+\int (\partial_t{{\overline{\rho}_{\varepsilon}}})^2\leq c\int \left\{{{\overline{\rho}_{\varepsilon}}}^2+(\nabla {{\overline{\rho}_{\varepsilon}}})^2+{{\overline{P}_{\varepsilon}}}^2\right\}+c
%\end{equation}
Equation for ${{\overline{P}_{\varepsilon}}}$: 
\begin{equation}\nonumber
\frac{\partial {{\overline{P}_{\varepsilon}}}}{\partial t}=\Delta {{\overline{P}_{\varepsilon}}} -{{\overline{P}_{\varepsilon}}}-\nabla {{\overline{\rho}_{\varepsilon}}}.
\end{equation}
%Energy estimate is (if multiply by ${{\overline{P}_{\varepsilon}}}$):
%\begin{equation}\label{en_2}
%\frac{\text{d}}{\text{d}t}\int {{\overline{P}_{\varepsilon}}}^2+\int |\nabla {{\overline{P}_{\varepsilon}}}|^2+\int {{\overline{P}_{\varepsilon}}}^2\leq \int |\nabla {{\overline{\rho}_{\varepsilon}}}|^2 
%\end{equation}
Energy estimate for ${{\overline{P}_{\varepsilon}}}$ (if multiply by ${{\overline{P}_{\varepsilon}}}|{{\overline{P}_{\varepsilon}}}|^2$):
\begin{equation}\nonumber
\frac{\text{d}}{\text{d}t}\int |{{\overline{P}_{\varepsilon}}}|^4 +12\int |\nabla {{\overline{P}_{\varepsilon}}}|^2|{{\overline{P}_{\varepsilon}}}|^2+4\int |{{\overline{P}_{\varepsilon}}}|^4=\int (\nabla {{\overline{\rho}_{\varepsilon}}}\cdot {{\overline{P}_{\varepsilon}}})|{{\overline{P}_{\varepsilon}}}|^2
\end{equation}
Estimate the right had side using integration by parts:
\begin{equation}\nonumber
\int (\nabla {{\overline{\rho}_{\varepsilon}}}\cdot {{\overline{P}_{\varepsilon}}})|{{\overline{P}_{\varepsilon}}}|^2\leq c\||\nabla {{\overline{P}_{\varepsilon}}}||{{\overline{P}_{\varepsilon}}}|\|_{L^2(\Omega)}\|{{\overline{\rho}_{\varepsilon}}}\|_{L^4(\Omega)}\|{{\overline{P}_{\varepsilon}}}\|_{L^4(\Omega)}.
\end{equation}
Thus, using 
\begin{equation}\nonumber
\||\nabla {{\overline{P}_{\varepsilon}}}||{{\overline{P}_{\varepsilon}}}|\|\|{{\overline{\rho}_{\varepsilon}}}\|_{L^4}\|{{\overline{P}_{\varepsilon}}}\|_{L^4}\leq \nu(\||\nabla {{\overline{P}_{\varepsilon}}}||{{\overline{P}_{\varepsilon}}}|\|^2+\frac{1}{c^2}\|{{\overline{\rho}_{\varepsilon}}}\|_{L^4}^4)+c\|{{\overline{P}_{\varepsilon}}}\|^4_{L^4} \end{equation}
\begin{equation}\nonumber
\leq \nu(\||\nabla {{\overline{P}_{\varepsilon}}}||{{\overline{P}_{\varepsilon}}}|\|^2+\frac{1}{c}\|\nabla {{\overline{\rho}_{\varepsilon}}}\|^2\|{{\overline{\rho}_{\varepsilon}}}\|^2)+c\|{{\overline{P}_{\varepsilon}}}\|^4\leq\nu(\||\nabla {{\overline{P}_{\varepsilon}}}||{{\overline{P}_{\varepsilon}}}|\|^2+\|\nabla {{\overline{\rho}_{\varepsilon}}}\|^2)+c\|{{\overline{P}_{\varepsilon}}}\|^4_{L^4} 
\end{equation}
we have 
\begin{equation}\label{en_3}
\frac{\text{d}}{\text{d}t}\int |{{\overline{P}_{\varepsilon}}}|^4 \leq \nu \|\nabla {{\overline{\rho}_{\varepsilon}}}\|^2+c\int|{{\overline{P}_{\varepsilon}}}|^4
\end{equation}
Adding \eqref{en_1} to \eqref{en_3} we get:
\begin{equation}\nonumber
\frac{\text{d}}{\text{d}t}\left[\int {{\overline{\rho}_{\varepsilon}}}^2 +\int {{\overline{P}_{\varepsilon}}}^4\right]\leq c\left[\int {{\overline{\rho}_{\varepsilon}}}^2 +\int {{\overline{P}_{\varepsilon}}}^4\right] 
\end{equation}
which proves uniqueness.

\subsection{Maximum princliple}
Consider the equation
\begin{equation}\label{as_original}
\partial_t \rho=\Delta \rho -\frac{W'(\rho)}{\ve^2}-P\cdot \nabla \rho +\lambda (t),\;\;\text{in } \Omega 
\end{equation}
with Neumann boundary conditions: 
\begin{equation}\nonumber
\nu \cdot \nabla \rho =0.
\end{equation}
Finctions $P$ and $\lambda$ are assumed to be given, $W'(\rho)=\frac{1}{2}\rho(1-\rho)(1-2\rho)$. 

\noindent{\bf Theorem.} Assume 
\begin{equation}\label{bound_on_ic}
0\leq \rho (x,0)\leq 1
\end{equation}
Then for all $t>0$
\begin{equation}\label{maximum_principle}
-2\ve^2\sup\limits_{\tau \in (0,t]}|\lambda(\tau)|\leq \rho(x,t)\leq 1+2\varepsilon^2\sup\limits_{\tau \in(0,t]}|\lambda (\tau)|.
\end{equation}

\medskip

\noindent{\it P r o o f:}\\
Denote $M:=\max\limits_{x\in \overline{\Omega},\tau\geq [0,t]}\rho(x,\tau)$ and assume that the maximum is attained in $x_0\in \Omega$ and $s_0 >0$. For such $x_0$ and $s_0$ we have: 
\begin{equation*}
\partial_t \rho\geq 0, \;\;\Delta\rho \leq 0,\;\;P\cdot \nabla \rho =0. 
\end{equation*}
Thus, 
\begin{equation*}
W'(M)\leq \ve^2\sup\limits_{s\in [0,t)}|\lambda (s)|.
\end{equation*}
Assume that $M>1$, then $W'(M)=\frac{1}{2}M(1-M)(1-2M)\geq \frac{1}{2}(M-1)$, so 
\begin{equation*}
M\leq 1+2 \ve^2\sup\limits_{s\in [0,t)}|\lambda (s)|.
\end{equation*}

\noindent Denote $m:=\min\limits_{x\in \overline{\Omega},s\geq [0,t]}\rho(x,s)$ and assume that the maximum is attained in $x_0\in \Omega$ and $s_0 >0$. For such $x_0$ and $s_0$ we have: 
\begin{equation*}
\partial_t \rho\leq 0, \;\;\Delta\rho \geq 0,\;\;P\cdot \nabla \rho =0. 
\end{equation*}
Thus, 
\begin{equation*}
W'(m)\geq -\ve^2\sup\limits_{s\in [0,t)}|\lambda (s)|.
\end{equation*}
Assume that $m<0$, then $W'(m)=\frac{1}{2}m(1-m)(1-2m)\leq \frac{1}{2}m$, so 
\begin{equation*}
-m\geq 2 \ve^2\sup\limits_{s\in [0,t)}|\lambda (s)|.
\end{equation*}

%\noindent The situation when $x\in \partial \Omega$ is considered in the similar manner using boundary conditions and the formula for the Laplacian on the boundary: 
%\begin{equation*}
%\Delta \rho|_{\partial \Omega}=\frac{\partial^2 \rho}{\partial \nu^2}+\frac{\partial^2 \rho }{\partial {\tau}^2}+\left(\frac{\partial \rho}{\partial \nu}\right)\nabla \cdot \nu.
%\end{equation*} 
%Here $\tau$ is tangent vector on $\partial \Omega$.  
\begin{acknowledgements}
We would like to thank referees for careful reading of the manuscript.
\end{acknowledgements}

\bibliographystyle{spmpsci}      % mathematics and physical sciences
\bibliography{cell}   % name your BibTeX data base

\end{document}